\documentclass{amsart} 
\usepackage{epsf} 
\usepackage{epstopdf} 
\usepackage{amssymb,latexsym, amsmath, amscd, array, graphicx,color} 

\usepackage{pdfsync} 
 
\usepackage{parskip}
\swapnumbers \numberwithin{equation}{section} 
\newtheorem{innercustomthm}{Theorem}
\newenvironment{customthm}[1]
  {\renewcommand\theinnercustomthm{#1}\innercustomthm}
  {\endinnercustomthm}
\newtheorem{innercustomcor}{Corollary}
\newenvironment{customcor}[1]
  {\renewcommand\theinnercustomcor{#1}\innercustomcor}
  {\endinnercustomcor}
\theoremstyle{plain} 
 
\newtheorem{thm}{Theorem}[section] 
\newtheorem{cor}[thm]{Corollary} 
\newtheorem{theorem}[thm]{Theorem} 
\newtheorem{lem}[thm]{Lemma} 
 
\newtheorem{prop}[thm]{Proposition}

\theoremstyle{definition} 
\newtheorem{defin}[thm]{Definition} 
 
\newtheorem{remark}[thm]{Remark}

\def\Im{\protect\operatorname{Im}}

\def\rot{\protect\operatorname{rot}}

\def\tb{\protect\operatorname{tb}}

\def\Z{{\mathbb Z}} 
 
\def\R{{\mathbb R}}

\def\1{\hbox{\rm\rlap {1}\hskip.03in{\rom I}}} 
\def\Bbbone{{\rm1\mathchoice{\kern-0.25em}{\kern-0.25em} {\kern-0.2em}{\kern-0.2em}I}} 

\def\pp{\medskip{\parindent 0pt \it Proof.\ }}

 \long
\def\forget#1\forgotten{} 

\begin{document}

\title[Loose Legendrian and Pseudo-Legendrian Knots]
{Loose Legendrian and Pseudo-Legendrian Knots in 3-Manifolds}
\author[P.~Cahn]{Patricia Cahn}
\address{Department of Mathematics,
Smith College,
Northampton, MA 01063, USA }
\email{pcahn@smith.edu}
\author[V.~Chernov]{Vladimir Chernov}
\address{Department of Mathematics, 6188 Kemeny Hall, Dartmouth College, Hanover, NH 03755, USA}
\email{Vladimir.Chernov@dartmouth.edu}

\begin{abstract}
We prove a complete classification theorem for loose Legendrian knots in an oriented 3-manifold, generalizing results of Dymara and Ding-Geiges.  Our approach is to classify knots in a $3$-manifold $M$ that are transverse to a nowhere-zero vector field $V$ up to the corresponding isotopy relation. Such knots are called $V$-transverse. A framed isotopy class is {\it simple\/} if any two $V$-transverse knots in that class which are homotopic through $V$-transverse immersions are $V$-transverse isotopic. We show that all knot types in $M$ are simple if any one of the following three conditions hold: $1.$ $M$ is closed, irreducible and atoroidal; or $2.$  the Euler class of the $2$-bundle $V^{\perp}$ orthogonal to $V$ is a torsion class, or $3.$ if $V$ is a coorienting vector field of a tight contact structure.  Finally, we construct examples of pairs of homotopic knot types such that one is simple and one is not. As a consequence of the $h$-principle for Legendrian immersions, we also construct knot types which are not Legendrian simple.

\end{abstract}

\maketitle
\section{Introduction}
We work in the smooth category. Throughout this paper $M$ is an oriented, connected $3$-manifold, which is not necessarily compact.  We fix an auxillary Riemannian metric on $M$.

Let $\xi$ be a cooriented contact structure on $M$.  For Legendrain knots in $(M,\xi)$ with well-defined rotation and Thurston-Bennequin numbers, one classical problem is: Given an ordered pair $(t,r)\in \mathbb{Z}\times\mathbb{Z}$, and a smooth knot type $\mathcal{K}$, classify the Legendrian knots $K\in \mathcal{K}$ such that $t=\text{tb}(K)$ and $r=\text{rot}(K)$.  This is sometimes referred to as the {\it botany problem} \cite{Etnyreovertwisted}.  The Thurston-Bennequin number $\text{tb}(L)$ of a Legendrian knot $L$ is defined when $L$ is zero-homologous in $M$.  The rotation number $\text{rot}(L)$ is defined either when $L$ is zero-homologous in $M$, or when $\xi$ is a trivializable 2-plane bundle.

We study the following {\it generalized botany problem}, which applies to all knot types in any contact manifold $M$ with a cooriented contact structure $\xi$: Given a connected component $\mathcal{FK}$ of the space of framed knots, and a connected component of the space of Legendrian immersed curves  $\mathcal{LC}$, classify the Legendrian knot types in $\mathcal{FK}\cap \mathcal{LC}$.  (The framing of a Legendrian knot $L$ is given by orthogonally projecting a coorienting vector field for the contact structure to the normal bundle of $L$.) This generalizes the botany problem because two smoothly isotopic Legendrain knots $L$ and $L'$ with well-defined Thurston-Bennequin numbers are isotopic as framed knots if and only if $\text{tb}(L)=\text{tb}(L')$, and similarly, two Legendrian knots with well-defined rotation numbers which are homotopic as immersed curves are homotopic as Legendrian immersed curves if and only if $\text{rot}(L)=\text{rot}(L')$.  We study this generalized botany problem by classifying knots transverse to a nowhere-zero vector field $V$ on $M$; we call these $V$-transverse knots.

Throughout the article, we use the following notation and terminology.  All spaces of knots and curves are equipped with the $C^\infty$-topology. Connected components of a space of knots or immersed curves are referred to as isotopy or homotopy classes, respectively.  A smooth isotopy class is denoted $\mathcal{K}$. An isotopy class of framed, Legendrian, or $V$-transverse knots is denoted $\mathcal{FK}$, $\mathcal{LK}$, or $\mathcal{VK}$, respectively.  A homotopy class of framed, Legendrian, or $V$-transverse immersed curves is denoted $\mathcal{FC}$, $\mathcal{LC}$, or $\mathcal{VC}$, respectively.  

 A Legendrian knot $L$ in an overtwisted contact manifold $(M, \xi)$ is  {\it loose} if it is contained in the complement of some overtwisted disk $D\subset M$.  Otherwise, $L$ is {\it non-loose}.  
 
  We give a complete solution to the generalized botany problem for loose Legendrian knots.
 
 Prior results of Dymara and Ding-Geiges suggest that, in the case of loose knots, classical invariants, and their generalizations discussed above, completely determine the Legendrian knot type. Dymara \cite{Dymara, DymaraOvertwisted} proved that if the 2-plane bundle $\xi$ is trivializable, and the smooth knot type $\mathcal{K}$ has infinitely many framings, then two loose Legendrian knots $L$ and $L'$ in $\mathcal{K}$ are isotopic as framed knots and homotopic through Legendrian immersions if and only if they are Legendrian isotopic.  The number of framings $|\mathcal{K}|$ of a smooth knot type $\mathcal{K}$ is the number of distinct isotopy classes of framed knots with underlying smooth knot type $\mathcal{K}$.  Ding and Geiges \cite{DingGeiges} generalized Dymara's theorem.  They proved that if $\mathcal{K}$ has infinitely many framings, and the connected component of the space of framed curves $\mathcal{FC}$ containing $L$ and $L'$ contains infinitely many distinct connected components of the space of Legendrian immersions (this is true, e.g., when $\xi$ is trivializable), then $L$ and $L'$ are isotopic as framed knots and homotopic through Legendrian immersions if and only if they are Legendrian isotopic.

  We first prove a best-possible generalization of the results of Dymara and Ding-Geiges.  Our generalization does not make any assumption on the number of framings of $\mathcal{K}$ or about the number of components of the space of Legendrian immersions in $\mathcal{FC}$.  Our proof uses an $h$-principle of Cieliebak and Eliashberg \cite{CieliebakEliashberg}.

\begin{customthm}{1}[cf. Theorems \ref{Vtransimpliesformal} and \ref{loosethm} in the text] \label{DymaraDingGeigesGeneralization} Let $M$ be a 3-manifold with a  cooriented contact structure $\xi$,  and let $D$ be an overtwisted disk in $M$. Let $L$ and $L'$ be two smoothly isotopic Legendrian knots in $M\setminus D$. 
Assume that the following three conditions hold:
\begin{enumerate}
\item $L$ and $L'$ are isotopic as framed knots,
\item $L$ and $L'$ are homotopic as Legendrian immersions,
\item $\Im \bar{i}_V=\Im \bar{h}_V,$ where $V$ is a coorienting vector field of $\xi$.

\end{enumerate}
Then $L$ and $L'$ are isotopic as Legendrian knots.\label{legendrian.thm}
\end{customthm}

In item (3), $\bar{i}_V$ and $\bar{h}_V$ are homomorphisms from the fundamental groups of the spaces of framed knots and immersed curves in $M$, respectively, to $\mathbb{Z}$; these homomorphisms are defined in Section ~\ref{kinkcancelling.sec}, and we compute them explicitly in many examples.  

Conversely, we construct examples of Legendrian knots which are {\it not} Legendrian isotopic, but are isotopic as framed knots and homotopic through Legendrian immersions.  These examples arise when $\Im \bar{i}_V\neq \Im \bar{h}_V$, and the knots in these examples can be chosen to be loose.

Our approach is to study a version of the botany problem for a more general class of knots, which we call $V$--transverse knots.    We say that a knot or immersed curve $K:S^1\rightarrow M$ is {\it $V$-transverse} if the velocity vector $\vec{K}'(t)$ and the vector $V_{K(t)}\in T_{K(t)}M$ span a $2$-plane for all $t\in S^1$.  In the special case where $V$ is a coorienting vector field for a contact structure $\xi$ on $M$, $K$ is called {\it pseudo-Legendrian}; this notion was introduced by Benedetti and Petronio \cite{BenedettiPetronio1, BenedettiPetronio2}.

\begin{figure}[htbp]
\includegraphics[width=7cm]{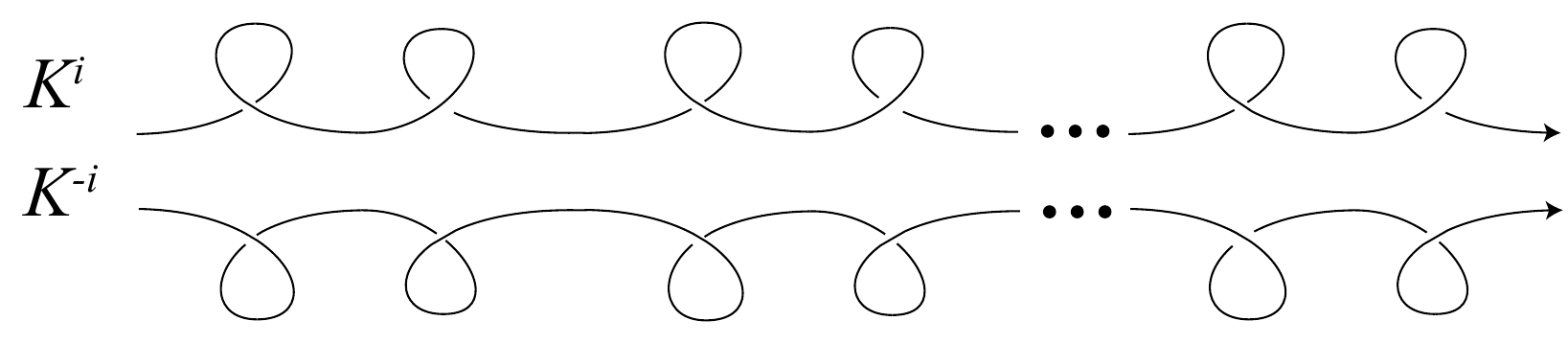}
\caption{The $V$-transverse knots $K^i$ and $K^{-i}$.  In this figure, $\frac{\partial}{\partial z}$ is pointing out of the page.}
\label{action.fig}
\end{figure}

We begin by proving Theorem ~\ref{customthm2}, a general classification theorem for $V$-transverse knots in oriented $3$-manifolds, which is the primary tool used in the proof of Theorem ~\ref{legendrian.thm}.  To state Theorem ~\ref{customthm2}, we introduce a local operation on $V$-transverse immersed curves.  Let $K$ be a $V$-transverse immersed curve in $M$.  Consider a coordinate chart $\phi:U\rightarrow \mathbb{R}^3$ such that $U\subset M$ contains an unknotted arc of $K$ and $V=\phi_*^{-1}\left(\frac{\partial}{\partial z}\right).$  Let $K^i$ denote the $V$-transverse knot obtained from $K$ by adding $i$ of the kink-pairs shown on the top in Figure \ref{action.fig} along this arc of $K$, if $i>0$, and $|i|$ of the kink-pairs on the bottom, if $i<0$.  In Lemma ~\ref{action.lem}, we prove that $i\mapsto K^i$ is a transitive $\mathbb{Z}$-action on the set of $V$-transverse homotopy classes in a given connected component of the space of framed curves in $M$, and also on the set of $V$-transverse isotopy classes in a given connected component of the space of framed knots in $M$.

\begin{figure}[htbp]
\includegraphics[width=9cm]{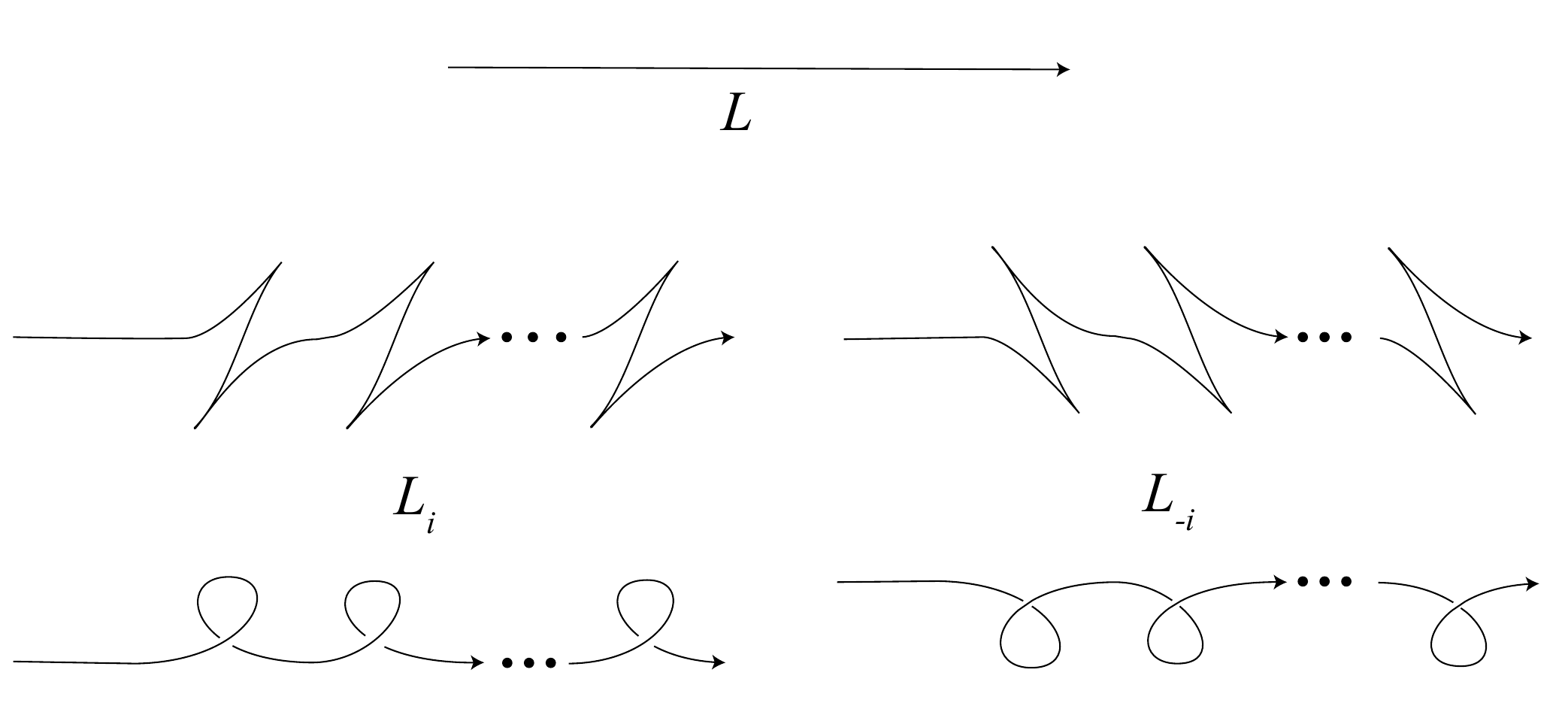}
\caption{The Legendrian knots $L_i$ and $L_{-i}$ in the front and Lagrangian projections.  In the Lagrangian projection, $\frac{\partial}{\partial z}$ is pointing out of the page. }
\label{zigzags.fig}
\end{figure}

This kinking operation is related to, but not the same as, the usual stabilization operation on Legendrian knots. The usual stabilization operation occurs in a chart on $(M,\xi)$ which is contactomorphic to $(\mathbb{R}^3,\xi_{\text{std}})$ where $\xi_{std}=\text{ker}(dz-ydx)$; such a chart exists at every point of $M$ by the Darboux Theorem.  Given a Legendrian knot $L$ in $(M,\xi)$ and positive integer $i$, let $L_i$ (respectively, $L_{-i}$) be the Legendrian knot obtained by performing $i$ positive (respectively, negative) stabilizations on $L$, as shown in Figure ~\ref{zigzags.fig}.   Working in the Lagrangian projection, and using the isotopy in Figure ~\ref{kink2.fig}, it is straightforward to check that $(L_{-i})^i$ and $L_i$ are isotopic as $V$-transverse knots, where again $V=\phi_*^{-1}\left(\frac{\partial}{\partial z}\right).$  Since $\frac{\partial}{\partial z}$ is a coorienting vector field for $\xi_{std}$, this $V$-transverse isotopy is a pseudo-Legendrian isotopy.

The above $\mathbb{Z}$-action appears in Theorem ~\ref{customthm2}, a classification theorem for $V$-transverse knots in a 3-manifold.  This theorem answers the $V$-transverse version of the botany problem: Classify the $V$-transverse knot types in a given framed isotopy class and a given connected component of the space of $V$-transverse immersed curves. (The framing of a $V$-transverse knot is given by orthogonally projecting $V$ to the normal bundle of $K$). 

 \begin{customthm}{2}[cf. Theorem \ref{classification.thm}\label{customthm2} in the text] Let $M$ be an oriented $3$-manifold and let $V$ be a non-vanishing vector field on $M$.  Let $K$ be a $V$-transverse knot in $M$, contained in the framed isotopy class $\mathcal{FK}$, the framed homotopy class $\mathcal{FC}$, and the $V$-transverse homotopy class $\mathcal{VC}$.  Let $\bar{h}_V:\pi_1(\mathcal{FC},K)\rightarrow \Z$ and $\bar{i}_V:\pi_1(\mathcal{FK},K)\rightarrow \Z$ be the homotopy and isotopy kink-cancelling homomorphisms discussed in detail in Section \ref{kinkcancelling.sec}.
 \begin{itemize}
 \item The set of $V$-transverse homotopy classes in $\mathcal{FC}$ is a $(\Z/\Im \bar{h}_V)$-torsor.
 \item The set of $V$-transverse isotopy classes in $\mathcal{FK}$ is a $(\Z/\Im \bar{i}_V)$-torsor.
 \item The set of $V$-transverse isotopy classes in $\mathcal{FK}\cap \mathcal{VC}$ is a $(\Im \bar{h}_V/\Im \bar{i}_V)$-torsor.
 \end{itemize}
 \end{customthm}

The maps $\bar{i}_V$ and $\bar{h}_V$ are defined as half the values of the Euler class of $V^\perp$ on the class of $H_2(M)$ realized by fundamental class of $S^1\times S^1$ pushed forward by the adjoint map $S^1\times S^1\to M$ of the loops in $\mathcal{FC}$ and $\mathcal{FK}$ respectively.  

The intuition behind Theorem ~\ref{customthm2} is as follows.  For a fixed $V$-transverse knot $K$, the knots $K^i$ defined by the kinking operation are always isotopic as framed knots.  Suppose $K$ is a $V$-transverse knot in $\mathbb{R}^3$ with $V=\frac{\partial}{\partial z}$.  In this case, the knots $K^i$ are all distinct as $V$-transverse immersed curves (and hence also as $V$-transverse knots), due to the fact that their projections to the $xy$-plane have different rotation numbers.  However, in an arbitrary oriented 3-manifold with nowhere-zero vector field $V$, it may happen that $K$ and $K^i$ are homotopic as $V$-transverse curves.  We prove this happens precisely when there is a framed self-homotopy $\alpha$ of $K$ along which the Euler class of the 2-plane bundle $V^\perp$ takes the value $2i$ (in this case $\bar{h}_V(\alpha)=i$).  Moreover, it may happen that $K$ and $K^i$ are distinct as $V$-transverse knots; we prove this happens precisely when there is no framed self-isotopy $\alpha$ of $K$ on which the Euler class of the 2-plane bundle $V^\perp$ takes the value $2i$ (in this case $i$ is not in the image of $\bar{i}_V$).  In particular, we get the following Corollary to Theorem ~\ref{customthm2}.

\begin{customcor}{3}[cf. Corollary ~\ref{numberofknottypes.cor} in the text]\label{customnumberofknottypes.cor} The number of $V$-transverse knot types in a given framed isotopy class which are homotopic through $V$-transverse immersions is the index of $\Im \bar{i}_V$ in $\Im \bar{h}_V$.
\end{customcor}

In $\mathbb{R}^3$ this index is 1.  In the examples in this paper, the index is either $1$ or infinite.  We do not know whether other values are possible.

We call a framed isotopy class $\mathcal{FK}$ {\it simple} if any two $V$-transverse knots $K$ and $K'$ in $\mathcal{FK}$ which are homotopic through $V$-transverse immersions are isotopic as $V$-transverse knots.   Theorem ~\ref{customthm2} allows us to determine exactly when a framed knot type is simple; see Corollary ~\ref{customequalimages.cor}.

\begin{customcor}{4}[cf. Corollary \ref{equalimages.cor} in the text]\label{customequalimages.cor}The framed isotopy class $\mathcal{FK}$ is simple if and only if $\Im \bar{h}_V=\Im \bar{i}_V$.  In particular, this is the case when $\Im \bar{h}_V=0$.
\end{customcor}

Next we describe three types of examples.  First, we describe manifolds $M$ and vector fields $V$ such that every framed knot type in $M$ is simple.

\begin{customthm}{5}[cf. Theorem \ref{simple.thm} in the text]  Let $V$ be a nowhere-zero vector field on an oriented $3$-manifold $M$ satisfying one of the following three conditions:
\begin{enumerate}
\item The Euler class $e_{V^\perp}\in H^2(M;\Z)$ is a torsion element, or in particular, if $e_{V^\perp}=0$.
\item The manifold $M$ is closed, irreducible and atoroidal.
\item $V$ is a coorienting vector field of a contact structure $\xi$ such that $(M,\xi)$ is tight, or more generally, such that $(M,\xi)$ is a covering of a tight contact manifold.
\end{enumerate}
Then every framed isotopy class in $M$ is simple.\label{simple.thm}
\end{customthm}

Second, we give examples of nonsimple classes. In our first set of examples, for an infinite family of vector fields $V_k$ with distinct Euler classes, we describe a framed homotopy class of immersed curves such that every framed knot type in this class is nonsimple.

\begin{customthm}{6}[cf. Theorem \ref{nonorientableex.thm} in the text] Let $M$ be an $S^1$-bundle over a non-orientable surface of genus $g\geq 1$ with oriented total space, and $\nu$ be the solid curve pictured in Figure \ref{example2.fig}. Let $\mathcal{FC}$ be any homotopy class (connected component of the space) of framed immersions in $M$ that contains a curve projecting to $\nu.$  Then for any nonzero $k\in \mathbb{Z}$ there exists a nowhere-zero vector field $V_k$ on $M$ such that no framed knot type in $\mathcal{FC}$ is simple.  In particular, for any $V_k$-transverse knot $K$ in $\mathcal{FC}$, the $V_k$-transverse knots $K$ and $K^k$ are homotopic through $V_k$-transverse immersions, isotopic as framed knots, and {\it not} isotopic through $V_k$-transverse knots.\label{nonsimple.thm}
\end{customthm}
As a corollary to Theorem ~\ref{nonsimple.thm}, we construct knot types which are not Legendrian simple.  Eliashberg ~\cite{Eliashbergovertwisted} proved that every 2-plane field is homotopic to an overtwisted contact structure.  Hence the vector field $V_k$ above can be chosen to be a coorienting vector field of an overtwisted contact structure $\xi_k$ on $M$.

\begin{customcor}{7} [cf. Corollary \ref{example1cor} in the text] Let $(M,\xi_k)$ be an oriented 3-manifold as in Theorem ~\ref{theoremexample1} with a cooriented contact structure $\xi_k$ homotopic to the 2-plane field $V_k^\perp$.   For each framed knot type $\mathcal{FK}\subset\mathcal{FC}$ containing a Legendrian representative $L$, the stabilized Legendrian knots $L_k$ and $L_{-k}$ are isotopic as framed knots, and homotopic as Legendrian immersions, but not isotopic as Legendrian knots; these Legendrian knots can be chosen to be loose.  
	\end{customcor}

Third, we describe a framed homotopy class of immersed curves which contains both simple and nonsimple framed knot types, again for an infinite family of vector fields $V_k$ with distinct Euler classes.  The proof of this theorem uses an invariant of properly immersed annuli in $4$-manifolds with boundary defined by Schneiderman ~\cite{Schneiderman}.

\begin{customthm}{8}[cf. Theorem \ref{orientableex.thm} in the text] Let $M$ be an $S^1$-bundle over an orientable surface of genus $g\geq 2$.  Let $\mathcal{FK}_1$ be the framed isotopy class of the $S^1$ fiber with any framing, and let $\mathcal{FK}_2$ be the framed isotopy class obtained from $\mathcal{FK}_1$ by a finger move around a curve projecting to the loop $l$ on $F$; see Figures \ref{simpleandnonsimple.fig} and \ref{K2.fig}. Then for any nonzero $k\in \Z$ there exists a nowhere-zero vector field $V_k$ such that the framed knot type of $\mathcal{FK}_1$ is simple in $(M, V_k)$ while the knot type $\mathcal{FK}_2$ is not.  In particular, there are $V_k$-transverse knots $K_2$ and $K_2^k$ in $\mathcal{FK}_2$ which  are homotopic through $V_k$-transverse immersions, isotopic as framed knots, and {\it not} isotopic through $V_k$-transverse knots.\label{simpleandnonsimple.thm}
\end{customthm}
Again, by choosing a contact structure $\xi_k$ on $M$ with coorienting vector field $V_k$, we construct examples of knot types which are not Legendrian simple.

\begin{customcor}{9} [cf. Corollary \ref{example2cor} in the text] Let $(M,\xi_k)$ be an oriented 3-manifold as in Theorem ~\ref{Theoremexample2}, with a cooriented contact structure $\xi_k$ homotopic to the 2-plane field $V_k^\perp$. Let $L$ be a Legendrian knot which is smoothly isotopic to the knot $K_2$ in Theorem ~\ref{Theoremexample2}. Then the stabilized Legendrian knots $L_k$ and $L_{-k}$ are isotopic as framed knots, and homotopic as Legendrian immersions, but not isotopic as Legendrian knots; these Legendrian knots can be chosen to be loose.  
\end{customcor}
The structure of the paper is as follows.  In Section \ref{classical.sec}, we review Trace's theorem for knots in $\mathbb{R}^3$ and explain why our classification theorem generalizes it.  In Section \ref{basic.sec}, we review basic facts of framed and $V$-transverse isotopy, and introduce the actions which appear in the main classification theorem.  In Section \ref{kinkcancelling.sec}, we introduce the homomorphisms $h_V$ and $i_V$ and prove the classification theorem.  In Section \ref{simple.sec}, we prove Theorem \ref{simple.thm}.  In Section \ref{specialloops.sec}, we discuss facts about the fundamental groups of the spaces of framed knots and immersions in $M$, which we use to construct the examples in the next two theorems.  Section \ref{vectorfields.sec} is a brief expository section which is helpful for visualizing the examples in the next two theorems.   Sections \ref{nonsimple.sec}, \ref{simpleandnonsimple.sec}, and \ref{legendrian.sec} contain proofs of Theorems \ref{nonsimple.thm}, \ref{simpleandnonsimple.thm}, and \ref{legendrian.thm}, respectively. The last section is an appendix on $h$-principles.

\section{Classical Invariants of $V$-Transverse Knots}\label{classical.sec}

Consider a knot $K$ in $\mathbb{R}^3$ transverse to the vertical vector field $V=\frac{\partial}{\partial z}$.  Consider the following two ``classical'' invariants of $K$: The rotation number of the projection of $K$ to the $xy$-plane, and the self-linking number of $K^V$, where $K^V$ is the knot $K$ framed by $V|_K$. Trace \cite{Trace} proved that such $V$-transverse knots are determined by their classical invariants.  We restate his theorem using our terminology.

\begin{thm}[Trace] Let $V=\frac{\partial}{\partial z}$. Two $V$-transverse knots $K$ and $L$ in $\mathbb{R}^3$ are isotopic as $V$-transverse knots if and only if
\begin{enumerate}
\item $K$ and $L$ are isotopic as smooth knots,
\item $K^V$ and $L^V$ have the same self-linking number, and
\item the projections of $K$ and $L$ to the $xy$-plane have the same rotation number.
\end{enumerate}\label{trace.thm}
\end{thm}

The second and third hypothesis above do not make sense in arbitrary $(M,V)$.  A rotation number of a $V$-transverse knot can only be defined given a trivialization of the $2$-plane bundle $V^\perp$, as the degree of the map which sends a point $t\in S^1$ to the normalized projection of the velocity vector $K'(t)$ to $V^\perp$.  The self-linking number is only defined for zero-homologous $K$.  When $K$ is not zero-homologous, one can instead use the affine self-linking invariant constructed by the second author \cite{ChernovFramed} which generalizes the ordinary self-linking number for zero-homologous knots and makes sense much more generally.

An even more general approach is to replace the first two hypotheses with the hypothesis that $K^V$ and $L^V$ be isotopic as framed knots, and replace the third hypothesis with the hypothesis that $K$ and $L$ are homotopic through $V$-transverse immersions.  These hypotheses make sense in any $(M,V)$ and are equivalent to Trace's hypothesis in $(\R^3,V=\partial/\partial z)$.

Thus a generalization of Trace's theorem should characterize when two $V$-transverse knots, which are homotopic through $V$-transverse immersions and isotopic as framed knots with framing given by $V$, are isotopic through $V$-transverse knots.

\section{Basic Properties of Framed and $V$-Transverse Isotopy}\label{basic.sec}

Let $M$ be a $3$-manifold.  A {\it framed curve} or {\it framed immersion} in $M$ is an immersion $C:S^1\rightarrow M$ together with a nonvanishing section of the normal bundle to $C(t)$ for each $t\in S^1$.  A {\it framed knot} is a framed curve which is also an embedding. The space of framed curves in $M$ always has two connected components corresponding to each connected component of the space of immersions of $S^1$ into $M$.  The space of framed knots may have finitely or infinitely many connected components corresponding to a given connected component of the space of knots in $M$ (a given framed isotopy class).  Both spaces are equipped with the $C^\infty$ topology.

\begin{prop}[Cf. \cite{CahnChernovSadykov}] Let $M$ be any $3$-manifold.  There are two components of the space of framed curves corresponding to each component of the space of unframed curves in $M$.
\end{prop}

\pp  Consider the Stiefel bundle of orthonormal $2$-frames $\xi$ over $M$.  The fiber is $V_2 \R^3=SO(3)=\R P^3$.  Let $C:S^1\rightarrow M$ be a framed curve in $M$.  Lift $C$ to a curve $\tilde{C}$ in $\xi$ using the frame given by $\{C'(t),v(t)\}$ where $v(t)$ is the framing vector of $C$ at time $t$.  Since each fiber of $\xi$ is an $\R P^3$ there is a canonical line bundle $E$ over $\xi$.  The first Stiefel-Whitney class of this bundle $w_1$ is an element of $H^1(\xi;\Z_2)$.  Its value on the lift of $C$ and the lift of $C$ with one extra twist of its framing differ by 1. Thus the number of connected components corresponding to each component of the space of unframed curves in $M$ is at least two.

Since there is an obvious homotopy between the curve with two extra twists of the framing and the original framed curve, the number of the connected components of the space of framed curves is at most two and hence actually equals two. 
\qed

Let $V$ be a nowhere-zero vector field on $M$.  A {\it $V$-transverse curve} or {\it $V$-transverse immersion} in $M$ is an immersion $C:S^1\rightarrow M$ such that $C'(t)$ and $V_{C(t)}$ span a $2$-plane for all $t\in S^1$.  A {\it $V$-transverse knot} is a $V$-transverse curve which is also an embedding.  Every $V$-transverse curve has a natural framing given by the orthogonal projection of $V_{C(t)}$ to the normal bundle of the curve at $C(t)$.

Recall from the introduction that, given a $V$-transverse knot $K$, there is a simple way to create a family $K^i$ of $V$-transverse knots via the kinking operation in Figure ~\ref{action.fig}.

\begin{figure}[htbp]\includegraphics[width=2in]{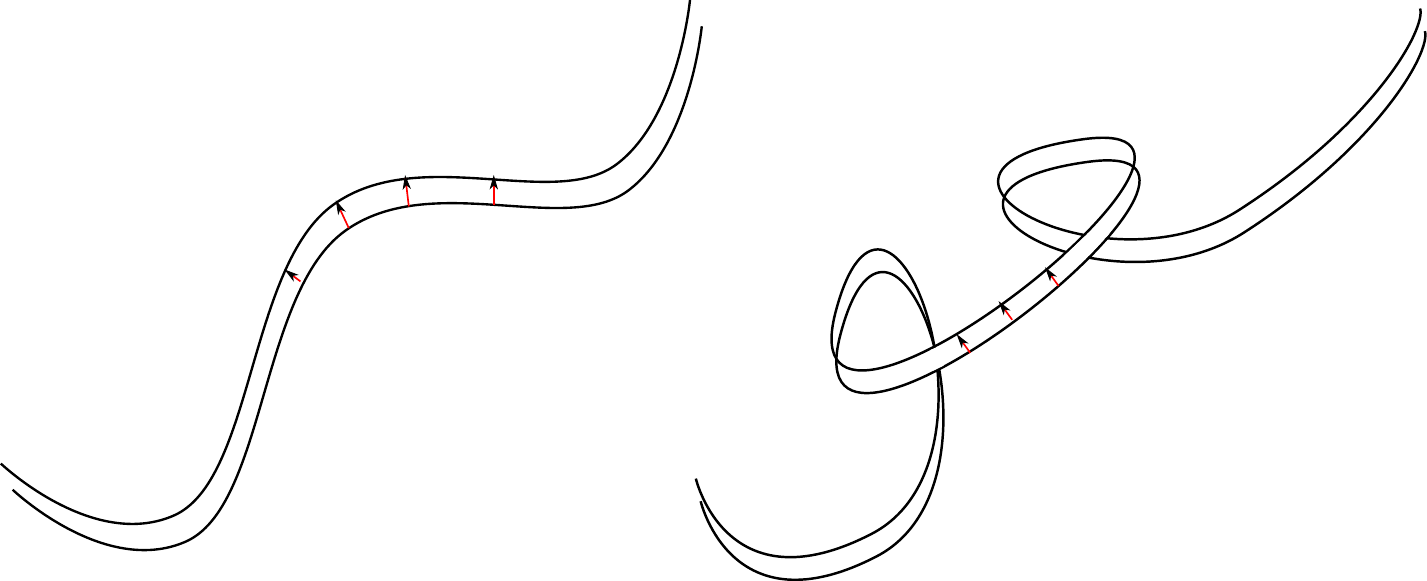}\caption{The $V$-transverse knots $K$ and $K^1$ are framed isotopic.}\label{framedisotopy.fig}
\end{figure}

\begin{prop}  The $V$-transverse knots $K^i$ are all isotopic as framed knots.
\end{prop}
\pp See Figure \ref{framedisotopy.fig}. \qed

\begin{prop}  Every framed isotopy class contains a $V$-transverse representative.
\end{prop}

\begin{figure}\includegraphics[width=3in]{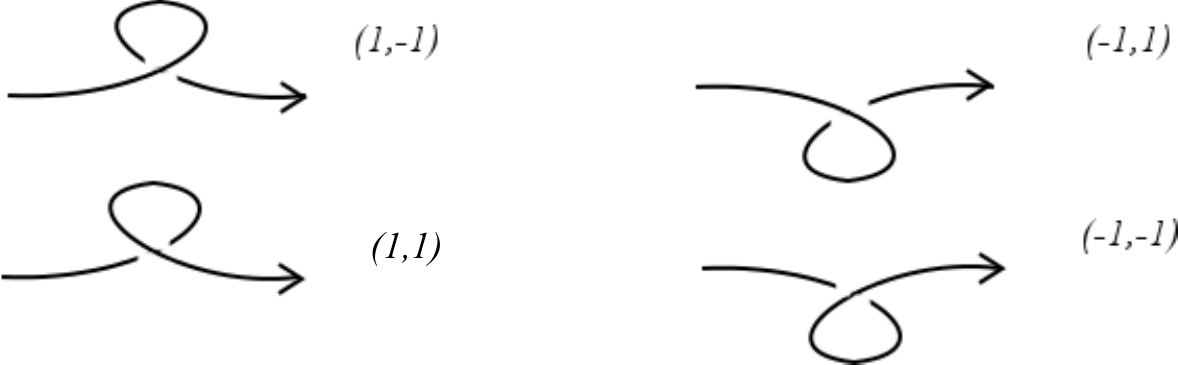}
\caption{Four different types of kinks with local contributions $(r,w)$ to the rotation number and writhe of the diagram.}
\label{fourkinks.fig}
\end{figure}

\pp Take an underlying unframed knot and make it transverse to a vector field $V$ so that $V$ gives some framing of the knot. Then add single kinks (see Figure ~\ref{fourkinks.fig}) until the corresponding framed knot has the desired framed knot type.\qed

Note that by Trace's theorem, the $K^i$ are all distinct as $V$-transverse knots, and even as $V$-transverse immersions, in $(\R^3,\partial/\partial z)$.  We will see that in other $(M, V)$ this need not be true.

\begin{lem}  Every framed isotopy can be $C^0$-approximated by a $V$-transverse isotopy.  In particular, if the $V$-transverse knots $K$ and $L$ are framed isotopic, then $L$ is $V$-transverse isotopic to $K^i$ for some integer $i$.\label{approximateisotopy.lem}
\end{lem}
\pp We can choose a set of coordinate charts $\{(U_i,\phi_i)\}_{i=1}^n$ for $M$ such that $V=\phi_{i*}^{-1}(\partial / \partial z)$ in each chart.  We will imitate the framed isotopy $K_t$ from $K=K_0$ to $L=K_1$ by a $V$-transverse isotopy $K^V_t$ in such a way that the knot $K^V_1$ agrees with $L$ outside some coordinate chart $(U_i,\phi_i)$, and inside that chart $L$ and $K^V_1$ differ by a collection of small kinks of four different types; see Figure ~\ref{fourkinks.fig}.  We will then argue that these kinks cancel via an isotopy in such a way that $K^V_1=L^i$.

In each chart we may assume that the projection of the framed isotopy $K_t$ to the $xy$-plane is a sequence of first, second, and third Reidemeister moves, in addition to ambient isotopy.

The second and third Reidemeister moves may appear in the projection of a $V$-transverse isotopy, but the first move does not appear, because the projection of a $V$-transverse isotopy to the $xy$-plane will always be an immersed curve.

There are four different kinds of kinks that may appear in a type 1 Reidemeister move, and these kinks are pictured in Figure~\ref{fourkinks.fig}.  Each kink is labeled by an ordered pair, where the first number is the contribution of the kink to the rotation number of the projection to the $xy$-plane, and the second is the local writhe number.

\begin{figure}[htbp]\includegraphics[width=4.5in]{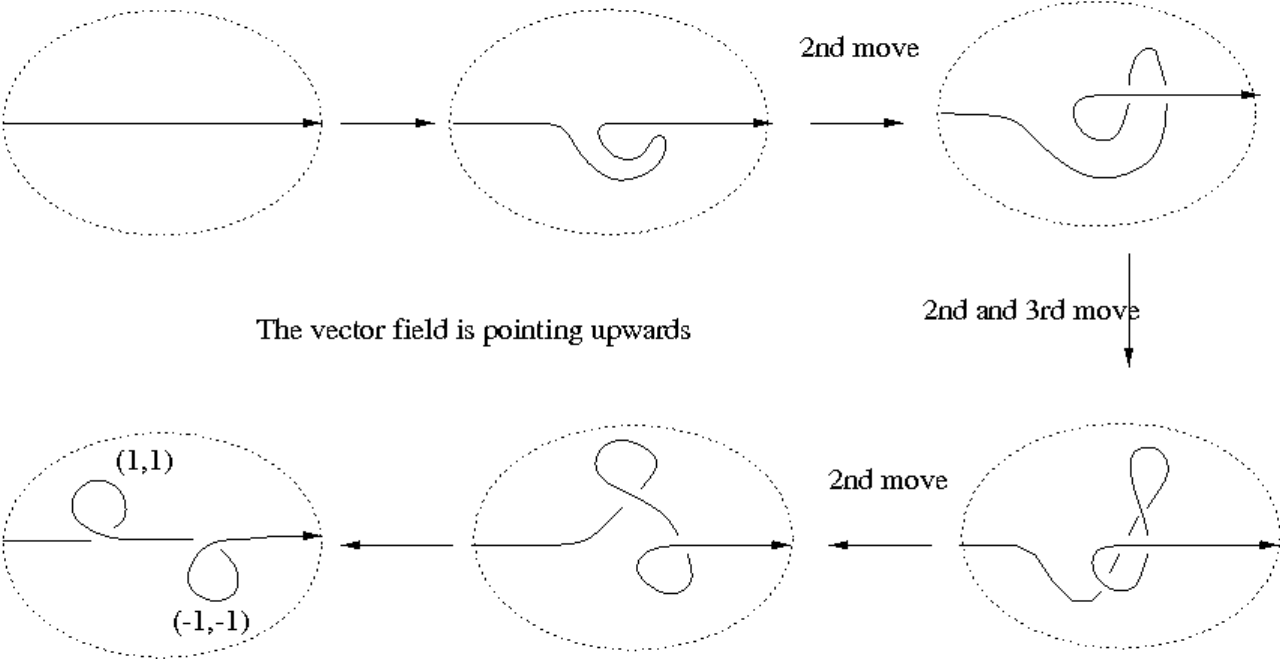}
\caption{Creation and cancelation of a pair of opposite kinks}
\label{kink2.fig}
\end{figure}

Pairs of kinks with opposite rotation number {\it and} opposite local writhe number can be created or cancelled by a $V$-transverse isotopy, see Figure~\ref{kink2.fig}.

Therefore if a type 1 move creates a kink of type $(\epsilon_1,\epsilon_2)$ during $K_t$, we instead create a pair of kinks $(\epsilon_1,\epsilon_2)$ and $(-\epsilon_1,-\epsilon_2)$ in $K^V_t$.  Then we make the extra kink of type $(-\epsilon_1,-\epsilon_2)$ very small and carry it along during the $V$-transverse isotopy.

If there is a type 1 move in $K_t$ which deletes a kink, we do not delete that kink in $K^V_t$ and instead make it small and carry it along during the $V$-transverse isotopy.

At the end of the isotopy $\overline{K}_t$ we see $L$ with many extra kinks.  We may slide these kinks along $L$ using a $V$-transverse isotopy so that they all appear in the same chart, and in an unknotted portion of $L$ in that chart. 

Let $a$ be the number of $(1,1)$ kinks, $b$ the number of $(-1,-1)$ kinks, $c$ the number of $(-1,1)$ kinks, and $d$ the number of $(1,-1)$ kinks.  Possibly by sliding kinks past one another, we cancel all pairs of kinks that have both opposite rotation number and opposite writhe. 

Hence we may assume that either $a$ or $b$ is equal to zero, and either $c$ or $d$ is equal to zero. For all $t$, the knots $K_t$ and $K^V_t$ are contained and isotopic in a thin solid torus $T_t$, whose core is $K_t,$ which we can identify with the standard solid torus in $\R^3$.   Since both $K_t$ and $K^V_t$ are framed isotopies, we can compare their self-linking numbers at each time $t$ after identifying $T_t$ with the standard solid torus in $\R^3$.  The difference between their self-linking numbers does not depend on the choice of identification of $T_t$ with the standard solid torus.  We call this number $s(t)$.  Since $s(0)=0$, we must have $s(1)=0$.  On the other hand each kink of type $(\epsilon_1, \epsilon_2)$ in $K^V_1$ contributes $\epsilon_2$ to the value of $s(1)$.  Hence $a+c=b+d$.

In the case where $a$ and $c$ are equal to $0$, we have $b+d=0$. But both $b$ and $d$ are nonnegative, so $a=b=c=d=0$.  In this case $K^V_1=L$.  This also occurs in the case where $b$ and $d$ are equal to $0$.

In the case where $a$ and $d$ equal $0$, we have $b=c$. In this case $K^V_1=L^{-b}$. In the case where $b$ and $c$ equal $0$, we have $a=d$ and $K^V_1=L^a$.  \qed

\begin{lem}  Every framed homotopy can be $C^0$-approximated by a $V$-transverse homotopy.  In particular, if the $V$-transverse curves $K$ and $L$ are framed homotopic then $L$ is $V$-transverse homotopic to $K^i$ for some integer $i$.\label{approximatehomotopy.lem}
\end{lem}

\pp Again, we cover $M$ with charts $(U_i,\phi_i)$ such that in each chart $V=\phi_{i*}^{-1}(\partial/\partial z)$.  Second and third Reidemeister moves, and crossing changes are $V$-transverse.  We adjust the first Reidemeister move as in the proof of Lemma~\ref{approximateisotopy.lem}.  At the end of our $V$-transverse homotopy, we are left with a copy of $L$ with extra kinks.  One can pass through a double point of a kink using a $V$-transverse homotopy, so we may cancel all pairs of kinks with opposite contributions to the rotation number, i.e., pairs of types $(\epsilon_1,\epsilon_2)$ and $(-\epsilon_1, \pm \epsilon_2)$.  We are left with kinks which all have the same local rotation number. Because $K$ and $L$ are in the same component of the space of framed curves the number of kinks remaining must be even. Now, we can pass through double points at vertices of the kinks to obtain $L^i$ for some $i\in \Z$.  
\qed

Let $V(\mathcal{F}K)$ denote the set of $V$-transverse knot types in the framed isotopy class $\mathcal{FK}$.  Let $V(\mathcal{FC})$ denote the set of $V$-transverse homotopy classes of $V$-transverse immersions in the framed homotopy class $\mathcal{FC}$.  We have now proven the following.

\begin{lem}\label{action.lem} The maps $\mathbb{Z}\times V(\mathcal{FC})\rightarrow V(\mathcal{FC})$ and $\mathbb{Z}\times V(\mathcal{FK})\rightarrow V(\mathcal{FK})$ defined by
$$i\cdot C\mapsto C^i$$
$$i\cdot K\mapsto K^i$$
define transitive actions on $V(\mathcal{FC})$ and $V(\mathcal{FK})$.
\end{lem}

The classification of $V$-transverse knots reduces to computing the stabilizers of these actions.
\section{The Kink-Cancelling Homomorphisms and the Classification Theorem}\label{kinkcancelling.sec}

Our goal is to measure the extent to which framed isotopic $V$-transverse knots can be homotopic as $V$-transverse immersions but not isotopic as $V$-transverse knots.  One can measure this using the Euler class of the $2$-plane bundle $V^\perp$.

In the following discussion we fix a framed isotopy class, or connected component $\mathcal{FK}$ of the space of framed knots and $\mathcal{K}$ the corresponding unframed isotopy class.  Let $\mathcal{FC}$ be the homotopy class of framed immersions (connected component of the space of framed curves), containing $\mathcal{FK}$, and $\mathcal{C}$ the corresponding homotopy class of unframed curves.

Suppose $a:S^1\rightarrow \mathcal{C}$ is a loop in the space of curves; that is, a self-homotopy of some curve $C$.  We regard $\alpha=[a]$ as an element of $\pi_1(\mathcal{C}, C)$.  Because $a(s)$ is a map $S^1\rightarrow M$, the map $a:S^1\rightarrow \mathcal{C}$ gives rise to a map of a torus, also called $a:S^1\times S^1\rightarrow M$, defined by $a(s,t)=a(s)(t)$.

Now we define a homomorphism $h_V:\pi_1(\mathcal{C},C)\rightarrow \Z$, which we call the {\it homotopy kink-cancelling homomorphism} by $$h_V(\alpha)=\frac{1}{2} e_{V^\perp}\left(a_*[S^1\times S^1]\right).$$

\begin{prop}  The map $h_V:\pi_1(\mathcal{C},C)\rightarrow \Z$ is a well-defined homomorphism.
\end{prop}
\pp  If $\alpha=[a_1]=[a_2]$ in $\pi_1(\mathcal{C},C)$ then $a_{1*}[S^1\times S^1]=a_{2*}[S^1\times S^1]$, so $h_V$ is well-defined and it clearly is a homomorphism. It is integer-valued because the value of $e_{V^\perp}$ is an even class, i.e.~it is $2\beta,$ for some $\beta\in H^2(M).$
\qed

Similarly, we define a homomorphism $i_V:\pi_1(\mathcal{K},K)\rightarrow \Z$, called the {\it isotopy kink cancelling-homomorphism}, by $$i_V(\alpha)=\frac{1}{2}e_{V^\perp}\left(a_*[S^1\times S^1]\right).$$

As before we have 
\begin{prop}  The map $i_V:\pi_1(\mathcal{K},K)\rightarrow \Z$ is a well-defined homomorphism.
\end{prop}

Note that in general, the map $inc_*:\pi_1(\mathcal{K},K)\rightarrow \pi_1(\mathcal{C},K)$ induced by the inclusion $\mathcal{K}\subset \mathcal{C}$ is neither one-to-one nor onto.  It is helpful to note that a loop in $\pi_1(\mathcal{C},K)$ ---i.e., a homotopy from $K$ to itself--- is in the image of $inc_*$ if it can be homotoped to an isotopy from $K$ to itself in the space of immersions of $S^1\rightarrow M$.  We will see an explicit example where $inc_*$ is not onto in Section~\ref{simpleandnonsimple.sec}.

The proof of the following proposition is straightforward.

\begin{prop} \label{imagecontainment.prop} Let $h_V$ and $i_V$ be the homotopy and isotopy kink-cancelling homomorphisms for the components $\mathcal{C}$ and $\mathcal{K}$ of the spaces of knots and curves, with basepoint $K$.  Then $\Im i_V\subseteq \Im h_V$.
\end{prop}

Now we give geometric interpretations of the homomorphisms above, and explain why we call them the {\it kink-cancelling} homomorphisms.  The proofs of the two lemmas below are similar, so we prove only the second.  

First we define two related homomorphisms, $\bar{h}_V$ and $\bar{i}_V$.  Given a loop $a\in \pi_1(\mathcal{FC},C)$, that is, a loop in space of {\it framed} curves, we can regard $a$ as an element of $\pi_1(\mathcal{C},C)$, a loop in the space of {\it unframed} curves, by forgetting the framing on $a$.  Then define $\bar{h}_V(a)=h_V(a)$.  The definition of $\bar{i}_V$ is analogous.  Note that the two lemmas below both involve the images of $\bar{h}_V$ and $\bar{i}_V$.  Also note that $\Im \bar{i}_V\subseteq \Im \bar{h}_V$ by the same argument as in Proposition ~\ref{imagecontainment.prop}.

\begin{lem}  Let $K$ be a $V$-transverse knot in $M$, and let $\mathcal{FC}$ be the homotopy class of framed curves containing $K$.  Let $\bar{h}_V$ be the corresponding kink-cancelling homomorphism defined on $\pi_1(\mathcal{FC},K)$.  Then $K$ and $K^k$ are homotopic as $V$-transverse immersions if and only if there is a framed self-homotopy $\alpha$ of $K$ such that $\bar{h}_V(\alpha)=k$. \label{Vtranshomo.lem}
\end{lem}

\begin{figure}[hbtp]
\includegraphics[width=.5in]{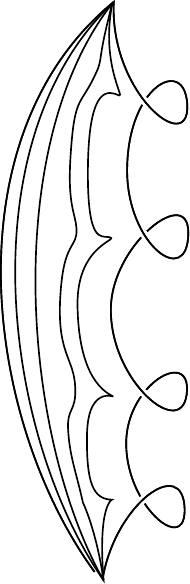}
\caption{An immersed disk in $\mathbb{R}^3$ formed by a framed isotopy from $K$ to $K^r$.}\label{kinkisotopy.fig}
\end{figure}

\begin{lem}  Let $K$ be a $V$-transverse knot in $M$, and let $\mathcal{FK}$ be the framed isotopy class of $K$.  Let $\bar{i}_V$ be the corresponding kink-cancelling homomorphism defined on $\pi_1(\mathcal{FK},K)$.  Then $K$ and $K^k$ are isotopic as $V$-transverse knots if and only if there is a framed self-isotopy $\alpha$ of $K$ such that $\bar{i}_V(\alpha)=k$.\label{Vtransiso.lem}
\end{lem}
\pp  {\it First we assume that $K(t)$ and $K^k(t)$ are isotopic as $V$-transverse knots, and show $k\in \Im \bar i_V$.\/}  Let $K^V_u$, with $u\in [0,1]$, be a $V$-transverse isotopy from $K$ to $K^k$, with $K^V_0=K$ and $K^V_1=K^k$.  Let $K^{\text{fr}}_v$, with $v\in [0,1]$, be the usual framed isotopy from $K^k$ to $K$ consisting of $k$ simultaneous applications of the move in Figure \ref{framedisotopy.fig}, with $K^\text{fr}_0=K^k$ and $K^\text{fr}_1=K$. We also assume that $K^\text{fr}_v(t)$ agrees with $K(t)$ for $t\in [0,1/2]$, so that the isotopy moves only an arc of $K$ as shown in Figure \ref{kinkisotopy.fig}.  Let $a:S^1\times S^1\rightarrow M$ be the self-isotopy of $K$ obtained by concatenating the two isotopies above:
$$a(s,t)=
\left\{\begin{array}{ll}K^V_{2s}(t)&\text{for }s\in [0,1/2]\\ K^\text{fr}_{2s-1}(t)&\text{for } s\in[1/2,1]\\ \end{array}\right.$$
There is a nowhere-zero section $\sigma$ of the $2$-plane bundle $V^\perp$ defined along the isotopy $K^V_u$, given by projecting the tangent vectors $(K^V_u)'(t)$ to $V^\perp$.  Consider the pullback $a^*(V^\perp)$ of the $2$-plane bundle to the torus $S^1\times S^1$.  We will show
$$e_{a^*(V^\perp)}([S^1\times S^1])=e_{V^\perp}(a_*[S^1\times S^1])=2k$$
by finding the obstruction to extending our section over the whole torus.
Choose a chart with $V=\partial/\partial z$, containing the part of the framed isotopy in which the $2k$ kinks are removed.  In this chart $\{\partial/\partial x,\partial/\partial y\}$ determine a trivialization of $V^\perp$, which is just a distribution of horizontal $2$-planes.  Pull the bundle $V^\perp$ back to the immersed disk $\phi:D^2=[1/2,1]\times [1/2,1] \rightarrow \mathbb{R}^3$ formed by the framed isotopy from $K$ to $K^r$, and equal to the image of $\alpha|_{[1/2,1]\times[1/2,1]}$ under the chart, shown in Figure \ref{kinkisotopy.fig}.  The degree of the map $\phi(\partial D^2)\mapsto \sigma$, computed with respect to the trivialization $\{\partial/\partial x,\partial/\partial y\}$, is $2k$.  Hence 
$$e_{a^*(V^\perp)}([S^1\times S^1])=e_{V^\perp}(a_*[S^1\times S^1])=2k$$
and $k\in \Im i_V$ as claimed.

Now assume $k\in \Im \bar{i}_V$.  We must construct a $V$-transverse isotopy from $K$ to $K^k$.  By assumption there is an element $\alpha=[a]\in \pi_1(\mathcal{FK}, K)$ such that $e_{V^\perp}(a_*[S^1\times S^1])=2k$.  View $a(s\times t)$ as a framed self-isotopy of $K$. By Lemma \ref{approximateisotopy.lem} there is a $C^0$-approximation of $a$ by a $V$-transverse isotopy $K_u^V$, such that $K^V_0=K$ and $K^V_1=K^i$ for some $i\in \Z$.  Concatenate this isotopy with the usual framed isotopy $K^{\text{fr}}_u$ from $K^i=K^\text{fr}_0$ to $K=K^\text{fr}_1$ given by $i$ applications of the move in Figure~\ref{framedisotopy.fig}, to get a map of the torus 
$$b(s,t)=
\left\{\begin{array}{ll}K^V_{2s}(t)&\text{for }s\in [0,1/2]\\ K^\text{fr}_{2s-1}(t)&\text{for } s\in[1/2,1]\\ \end{array}\right.$$
The maps $a$ and $b:S^1\times S^1\rightarrow M$ are $C^0$-close, so $a_*[S^1\times S^1]=b_*[S^1\times S^1]\in H_2(M,\Z)$. Thus $e_{V^\perp}(b_*[S^1\times S^1])=2i$, so $i=k$, and we have a $V$-transverse isotopy from $K$ to $K^i$.
\qed

The following theorem gives a complete classification of $V$-transverse knots.

Recall that $V(\mathcal{FK})$ is the set of $V$-transverse knot types in the framed isotopy class $\mathcal{FK}$, and $V(\mathcal{FC})$ is the set of $V$-transverse homotopy classes in the framed homotopy class $\mathcal{FC}$.  Let $V(\mathcal{FK}\cap \mathcal{VC})$ denote the set of $V$-transverse knot types in $\mathcal{FK}\cap \mathcal{VC}$.

Theorem \ref{classification.thm} below is Theorem~\ref{customthm2} from the Introduction.

\begin{thm}  Let $K$ be a $V$-transverse knot in $M$, contained in the framed isotopy class $\mathcal{FK}$, the framed homotopy class $\mathcal{FC}$, and the homotopy class of $V$-transverse immersed curves $\mathcal{VC}$.  Let $\bar{h}_V:\pi_1(\mathcal{FC},K)\to \Z$ and $\bar{i}_V:\pi_1(\mathcal{FK},K)\to \Z$ be the homotopy and isotopy kink-cancelling homomorphisms  for framed curves and knots respectively. Then
\begin{itemize}
\item $V(\mathcal{FC})$ is a $(\Z/\Im \bar{h}_V)$-torsor;
\item $V(\mathcal{FK})$ is a $(\Z/\Im \bar{i}_V)$-torsor;
\item $V(\mathcal{FK}\cap \mathcal{VC})$ is a $(\Im \bar{h}_V/\Im \bar{i}_V)$-torsor.
\end{itemize}\label{classification.thm}
\end{thm}

Of particular interest are the following corollaries, which illustrate how the classification may differ from the familiar example of $(\mathbb{R}^3,\partial/\partial z)$. They are Corollaries
~\ref{customnumberofknottypes.cor} and \ref{customequalimages.cor} from the Introduction.

\begin{cor}\label{numberofknottypes.cor}

The number $|V(\mathcal{FK}\cap \mathcal{VC})|$ of $V$-transverse knot types in a given framed isotopy class which are homotopic through $V$-transverse immersions is the index of $\Im \bar{i}_V$ in $\Im \bar{h}_V$.

\end{cor}

In $(\mathbb{R}^3,\partial/\partial z)$ this index is 1.  In the examples in this paper, the index is either $1$ or infinite.  We do not know whether other values are possible.

\begin{cor}\label{equalimages.cor}  
The framed isotopy class $\mathcal{FK}$ is simple if and only if $\Im \bar{h}_V=\Im \bar{i}_V$.  In particular, this is the case when $\Im \bar{h}_V=0$.
\end{cor}

\section{Some very basic examples}\label{surfaces.sec}
Let $F$ be an oriented surface and $M=F\times \mathbb{R}$.  Let $V=\partial/\partial z$ where $z$ is the $\mathbb{R}$ coordinate.  As in the case $F=\mathbb{R}^2$ (the setting of Trace's theorem), $V$-transverse knots are described by regular knot diagrams on $F$ up to the second and third Reidemeister moves, and ambient isotopy.  The framing is the blackboard framing given by the $\mathbb{R}$ factor.  We consider closed surfaces in these examples, since in those cases, $V^\perp$ (which can be identified with $TF\times\{0\}$) is not necessarily trivializable, so our theory is more interesting.

\begin{figure}[htbp]\includegraphics[width=4in]{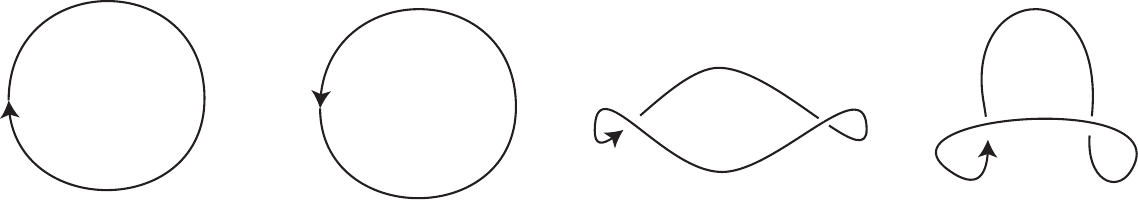}
\caption{An isotopy fro $K$ to $K^2$ on $S^2$.}
\label{sphereisotopy.fig}
\end{figure}

{\bf Example 1.  $F=S^2$.} There are two framed isotopy classes of framed
knots corresponding to each isotopy class of unframed knots.  Let $K$
be a small circle embedded in $S^2$ with framing given by $V$.  We will show that the framed knot type $K$ is simple.  There
is a framed isotopy taking $K$ to $K^2$, pictured in Figure
~\ref{sphereisotopy.fig}.  Note that it is important to keep in mind $K^2$ means $K$ with two extra pairs of kinks, as defined by our action, not $K$ with two extra twists of its framing.  The first step takes a small loop and
isotopes it around the back of the sphere to reverse its orientation.
The next steps are a Reidemeister 2 move, followed by an ambient
isotopy.  The torus swept out by this isotopy followed by the usual
framed isotopy from $K^2$ back to $K$ is homologous to the
fundamental class $[S^2]$.  Furthermore $TM=TS^2\times \mathbb{R}$ and
$V^\perp$ can be identified with the bundle $TS^2\times \{0\}$.  We
see that there is a self-isotopy $\alpha$ of $K$ such that
$e_{V^\perp}(\alpha_*[S^1\times S^1])=\chi(S^2)=2$.  Thus $\bar{i}_V$ is onto, and $\Im
\bar{i}_V =\mathbb{Z}=\Im \bar{h}_V$.  By Corollary
\ref{equalimages.cor} two knots in the framed isotopy class of $K$ are
isotopic through $V$-transverse knots if and only if they are
homotopic through $V$-transverse immersions.  The same can be said for
the framed isotopy corresponding to $K$ with an extra twist of its framing.

{\bf Example 2.  $F=T^2$.}  In this case all classes are simple because $V^\perp$ is trivial.  Note that there are still interesting self-isotopies of knots on the torus (i.e., self-isotopies such that the corresponding torus is not zero-homologous).  For example, if $K$ is the meridian, there is a self-isotopy of $K$ which sweeps out the whole torus.

{\bf Example 3.  $F$ is a closed surface of genus $g\geq 2$.}  So $F$ is a $K(\pi_1)$. The only abelian subrgoups of $\pi_1(F)$ are infinte cylic by the Preissman's Theorem~\cite{DoCarmo}, and $\pi_1(S^1\times S^1)=\Z\oplus \Z$ is abelian. Hence a map $S^1\times S^1\to F$ factors through $S^1$ and for a map $S^1\times S^1\to F\times\R$ the image of the fundamental class of the torus is $0 \in H_2(F\times \R)=H_2(F)$. Thus the homomorphisms $\bar i_V$ and $\bar h_V$ are zero homomorphisms.

\section{Simple Knot Types}\label{simple.sec}

In this section we point out interesting properties of $M$ and $V$ which cause {\it all} framed isotopy classes of knots in $M$ to be simple automatically.

A $3$-manifold $M$ is {\it irreducible} if every $2$-sphere embedded in $M$ bounds a ball.  $M$ is {\it atoroidal} if there are no maps $f:S^1\times S^1\rightarrow M$ which are $\pi_1$-injective.

\begin{theorem}  Assume that the pair $(M,V)$ satisfies any one of the following conditions:
\begin{enumerate}
\item The Euler class $e_{V^\perp}\in H^2(M;\Z)$ is a torsion element, or in particular, if $e_{V^\perp}=0$;
\item The manifold $M$ is closed, irreducible and atoroidal, which includes the case where $M$ is equipped with a Riemannian metric of negative sectional curvature;
\item $V$ is a coorienting vector field of a contact structure $\xi$ such that $(M,\xi)$ is tight, or more generally, such that $(M,\xi)$ is a covering of a tight contact manifold.
\end{enumerate}
Then every framed isotopy class in $M$ is simple.\label{simple.thm}
\end{theorem}
\pp  It suffices to show that if any one of these three conditions holds, then for every $\alpha:S^1\times S^1\rightarrow M$ we have $e_{V^\perp}(\alpha_*[S^1\times S^1])=0$ (implying both kink-cancelling homomorphisms are zero).  If condition $(1)$ holds this is certainly true.  If condition $(2)$ holds this is true because $\alpha_*[S^1\times S^1]=0$; see for example \cite[p. 2784]{ChernovTransverse}.  If condition $(3)$ holds then $e_{V^\perp}=e_{\xi}$ and the desired statement was proven in \cite[Corollary 3.10]{ChernovTransverse}.  (Note that if $\alpha:S^1\times S^1\rightarrow M$ is an embedding, then $e_\xi(\alpha_*[S^1\times S^1])=0$ by a Bennequin type inequality of Eliashberg \cite[Theorem 2.2.1]{Eliashberg}.\qed

\section{Some special loops in the space of framed immersions of $S^1\rightarrow M^3$}\label{specialloops.sec}

In this section we recall three important elements of $\pi_1(\mathcal{FC},K)$ and $\pi_1(\mathcal{C},K)$ from \cite{ChernovFramed}.  These loops will be used to construct examples of $V$-transverse knots which are framed isotopic, homotopic through immersions, and not $V$-transverse isotopic. 

\subsection{The number of framings of a knot.} \label{numberofframings.sec} First we recall a result about framed knots. It follows from the existence of the self-linking number that the number $|K|$ of framed knots in $S^3$ with given underlying zero-homologous knot $K$ is infinite. The second author previously defined affine self-linking invariants and used them~\cite[Theorem 2.4]{ChernovFramed} to show that $|K|$ is infinite for every knot in an orientable manifold unless the manifold contains a connected sum factor of $S^1\times S^2$. The knot $K$  need not be zero-homologous and the manifold is not required to be compact. In our work with Sadykov~\cite{CahnChernovSadykov} we used the results of McCullough~\cite{McCullough} and strengthened the above result.
We showed that 

\begin{lem}\label{finitelymanyframings} Let $M$ be a not necessarily compact orientable 3-manifold.  Given a knot $K$ in $M$ we have $|K|=\infty$ unless $K$ intersects a nonseparating $2$-sphere at exactly one point  in which case $|K|=2.$
\end{lem} 

Note that if $K$ intersects a nonseparating sphere at exactly one point then $M$ contains $S^1\times S^2$ as a connected sum factor.

\subsection{The framing loop $\gamma_{\text{fr}}$}
\begin{figure}[htbp]
 \begin{center}

\includegraphics[width=10cm]{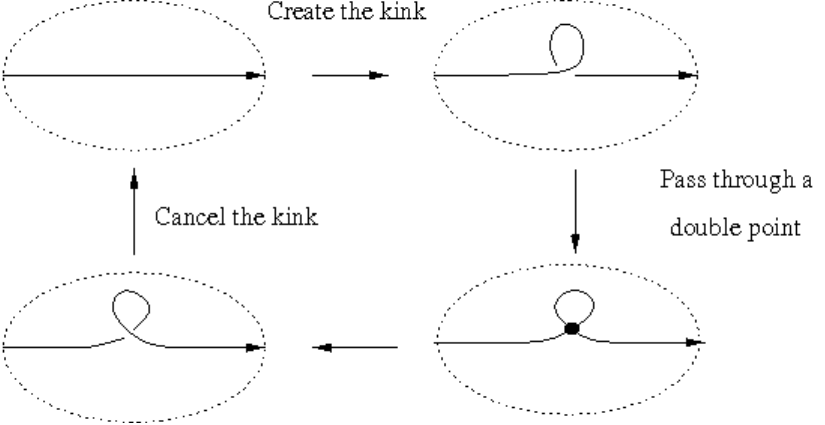}  
\end{center}
\caption{The loop $\gamma_{\text{fr}}$.} 
\label{obstruction.fig}
\end{figure}

Let $\gamma_{\text{fr}}$ be the element of $\pi_1(\mathcal{C}, K)$ pictured in Figure \ref{obstruction.fig}.  We call it the {\it framing loop} because it locally changes the framing of a framed knot by adding two full twists.  We say locally because Lemma \ref{finitelymanyframings} implies that in some cases the resulting knot is actually framed isotopic to $K$.  The loop $\gamma_{\text{fr}}$ is {\it not} an element of $\pi_1(\mathcal{FC},K)$; rather, it is a path from $K$ to $K^2$, the framed knot $K$ with two extra twists added to its framing.  $K$ and $K^2$, as stated above, may or may not be isotopic as framed knots, but nevertheless they are still always distinct points in $\mathcal{FC}$.

\subsection{The rotation loop $\gamma_{\text{rot}}$}
Let $\gamma_{\text{rot}}$ be the element of $\pi_1(\mathcal{FC}, K)$ induced by one full rotation of the parameterizing circle of $K$.

\subsection{The fiber loop $\gamma_{\text{fib}}$}
This loop is defined in the special case where $p:M^3\to F$ is an $S^1$-bundle over a (not necessarily orientable) surface $F$ and $p(K)$ is an orientation-preserving loop on $F.$ Since $M$ is oriented and $p(K)$ is orientation-preserving, we can orient the $S^1$-fibers containing the points of $K$ so that this orientation continuously depends on the point of $K$. Note that 
if a double point of $p(K)$ separates $p(K)$ into two orientation reversing loops then the corresponding two points of $K$ give different orientations of the $S^1$-fiber. 

 Let $\gamma_{\text{fib}}$ be the homotopy of $K$ that slides every point $K(t)$ of $K$ around the fiber that contains $K(t)$ with unit velocity, in the direction specified by the orientation of the fiber.  The homotopy $\gamma_{\text{fib}}$ is an element of $\pi_1(\mathcal{FC},K)$, but may or may not be an element of $\pi_1(\mathcal{K},K)$.  

\subsection{The loops $\gamma_{\rho}$}\label{gammarho}
This loop is defined only in the special case where $p: M\to F$ is an $S^1$-bundle over a (not necessarily orientable) surface, $\rho=[r]\in \pi_1(F)$ is the class of an orientation-preserving loop $r$ on $F$ based at $p(K(1))$, and $K$ is an oriented $S^1$-fiber.  We consider the framed isotopy $\gamma_{\rho}$ such that at each time moment $t$ the underlying loop $\gamma_{\rho}(t)$ is the $S^1$-fiber over $p(r(t)).$
There are many choices for such an isotopy, with two such differing by a power of $\gamma_{\text{fib}}$.  For our purposes this is not a problem.  The isotopy $\gamma_{\rho}$ may be regarded as an element  of $\pi_1(\mathcal{C},K)$, or $\pi_1(\mathcal{K},K)$.

\section{Visualizing vector fields with a given Euler class of $V^{\perp}$}\label{vectorfields.sec}

This section is a review of the Pontryagin-Thom construction and is expository, but useful for explicitly constructing and visualizing the vector fields which appear in the examples throughout the rest of the paper. We follow Geiges~\cite[Section 4.2]{Geiges}.

Our goal will be to construct and visualize a vector field for which the Euler class $e_{V^\perp}$ is $2k[\tilde{d}]$ for some link $\tilde{d}$ in $M$.

It is easier to begin with the reverse direction, and explain how to
construct and visualize the Poincar\'e dual of the Euler class of $V^{\perp}$.  Fix a trivialization of $TM\simeq M\times \mathbb{R}^3$.  Let
$V$ be a nowhere-zero vector field on $M$.  Then using our fixed
trivialization, $V$ corresponds to a map $f_V:M\rightarrow S^2$.
Choose a basis for $\mathbb{R}^3$ such that the north pole $N$ is a
regular value for $f_V$.  Let $\tilde{d}=f_V^{-1}(N)$, which is a link
in $M$.  Then the Euler class $e_{V^\perp}$ is Poincar\'e dual to $2
[\tilde{d}].$  To see why, recall that the Euler class of $V^\perp$ is
Poincar\'e dual to the zero set of a generic section of $V^\perp$.  To
get a generic section of $V^\perp$, project $V$ to the $xy$-plane
(possibly after a homotopy of $V$ to ensure it is generic, and to
ensure that the south pole $S$ is also a regular value of $f_V$),
according to our trivialization. The resulting section is zero along
$f^{-1}(N)\cup f^{-1}(S)$, the preimages of the poles.  The preimage
of the north pole is $\tilde{d}$, and one can check that the preimage
of the south pole is homologous to the preimage of $N$.  Thus
$2[\tilde{d}]$ is Poincar\'e dual to $e_{V^\perp}$.

Now we consider the question we are interested in: how to construct a
vector field such that $e_{V^\perp}$ is Poincar\'e dual to
$2k[\tilde{d}]$.  We start with the case $k=1$.  Again fix a
trivialization of $TM$.  Define $V$ to be $N$ along $\tilde{d}$.  Fix
a tubular neighborhood $T=\tilde{d}(t)\times D^2$ of $\tilde{d}$.  Define $V$ to be
$S$ along $M-\mathring{T}$.  Now define $V$ along the rest of $T$ by
mapping each cross-section $\tilde{d}(t)\times \mathring{D}^2$ to the open disk $S^2-\{S\}$ such that
the center of the disk maps to $N$ (as already noted).  Perturb this
vector field slightly so that $S$ is a regular value of $f_V$ (for
example, take $V$ to be a point close to $S$ in $M-\mathring{T}$).
Now $f_V^{-1}(S)$ runs parallel to $\tilde{d}$ in $T$, and the Euler
class $e_{V^\perp}$ is $2[\tilde{d}]$.

To construct $V_k$ dual to $2k[\tilde{d}]$, repeat the above process
with $k$ parallel copies of $\tilde{d}$.

\section{Examples of nonsimple Knot Types}\label{nonsimple.sec}
Given $k\in\Z$, in every $S^1$-bundle $M$ over a non-orientable surface $F$ of genus at least $1$ with an oriented total space, we can construct a nowhere-zero vector field $V=V_k$ on $M$, and a $V$-transverse knot $K$, such that $K$ and $K^k$ are framed isotopic, homotopic as $V$-transverse immersions, and {\it not} isotopic as $V_k$-transverse knots.

For such a knot $K$, since $K$ is isotopic to  $K^k$ through $V$-transverse knots, it must be the case that $k\in \Im \bar{h}_V$.  And, since $K$ and $ L$ are not homotopic through $V$-transverse immersions, we must have that $k\notin \Im \bar{i}_V$.

In particular, we seek a knot $K$ and a vector field $V$, such that there exists a framed homotopy from $K$ to itself, such that the pullback of the bundle $V^\perp$ to the corresponding torus has Euler number $2k$, {\bf but} for {\bf every} framed isotopy from $K$ to itself, the pullback of the bundle $V^\perp$ to the corresponding torus has Euler number $0$.
\begin{figure}[htbp]\includegraphics[width=3in]{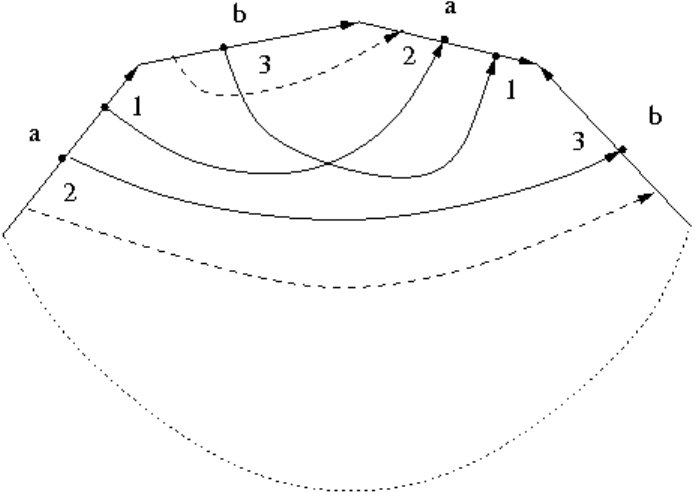}
\caption{$V$ is the vector field on $M$ such that the Euler class of the $2$-plane bundle $V^\perp$ is Poincar\'e dual to a lift $\tilde{d}$ of the dotted curve $d$ to $M$. The knot $K$ is a lift of the solid curve $\nu$.}\label{example2.fig}
\end{figure}

For any fixed $k\in \mathbb{Z}$, let $V=V_k$ be any vector field on $M$ such that the Euler class of the $2$-plane bundle $V^\perp$ is Poincar\'e dual to $2k[\tilde{d}]$, where $\tilde{d}$ is a curve in $M$ that projects to  the dotted curve $d$ in Figure \ref{example2.fig}.  Let $K$ be any knot in $M$ that projects to the solid curve $\nu$ in Figure \ref{example2.fig}. We assume that $K$ is perturbed slightly so it is $V$-transverse.  This does not determine a unique $V$-transverse isotopy class, but that is not a problem.  Let $\mathcal{FC}$ be the connected component of the space of framed immersions containing $K$. Let $K'$ be any other $V$-transverse knot in $\mathcal{FC}$.

\begin{lem}[\cite{ChernovFramed} proof of Lemma 6.11]  Let $\alpha\in \pi_1(\mathcal{C},K)$.  Then $\alpha$ contains a representative of the form $\gamma_{\text{fib}}^i\gamma_{\text{rot}}^j\gamma_{\text{fr}}^k$.
\label{example2standardform.lem}
\end{lem}

\pp Pick a representative $a\in \alpha$.  Since $a$ is a loop in the space of immersions, we can lift $a$ to a loop also called $a\in \pi_1(STM,\bar{K})$ where $\bar{K}$ is the lift of $K$ to the spherical tangent bundle $STM$.  Now we regard $a$ as a map $S^1\times S^1\rightarrow STM$ and consider the possible values of $\tau(\alpha)$, the trace of the basepoint of $a$ in $\pi_1(M, K(1))$.  For the precise definition of $\tau$ see Section \ref{hprinciples.sec}. The image of $\tau:\pi_1(\Omega STM,\bar{K})\rightarrow \pi_1(STM, \bar{K}(1))=\pi_1(M, K(1))$ is the centralizer of $K$. From the long exact sequence of the bundle $S^1\hookrightarrow M \rightarrow F$, and from the fact that the centralizer of $\nu\in \pi_1(F)$ is the cyclic subgroup generated by $\nu$, we conclude that any element of the centralizer of $K$ in $\pi_1(M, K(1))$ is of the form $f^iK^j$ where $f$ is the class of the $S^1$-fiber.  On the other hand, $\tau(\gamma_\text{fib})=f$ and $\tau(\gamma_\text{rot})=K$, so by Proposition \ref{sametrace.prop}, $a$ is homotopic to a representative of the form $\gamma_{\text{fib}}^i\gamma_{\text{rot}}^j\gamma_{\text{fr}}^k$. \qed

For the next Lemma, we introduce an invariant $\delta$ of knots in $\mathcal{C}$.  Let $K_s$ be a singular knot with one double point at $s\in M$ and view $K_s$ as a pair of maps $(K_{s,1},K_{s,2})\in \pi_1(M,s)\times \pi_1(M,s)$.  Let $\sigma(K_s)=1$ if both $K_{s,i}$ are noncontratible and let $\sigma(K_s)=0$ otherwise.  We say a loop or path $\gamma:[0,1]\rightarrow \mathcal{C}$ is {\it generic} if, whenever $\gamma(t)$ is a singular knot, $\gamma(t)$ has exactly one transverse double point and no other multiple points, and the set of times such that $\gamma(t)$ is singular is a discrete set $\{t_1,\dots,t_n\}$.  In particular any loop in $\pi_1(\mathcal{C},K)$ has a generic representative.  

A transverse double point $s$ of a singular knot can be resolved in two different ways. We say that a resolution of a double point is
positive (resp. negative) if the tangent vector to the
first strand, the tangent vector to the second strand, and the vector from
the second strand to the first form a positive $3$-frame. This does 
not depend on the order of the strands.

We assign a sign to each singular knot $\gamma(t_i)$ as follows: if, for $t_i<t^+<t_{i+1}$, $\gamma(t^+)$ is obtained from the singular knot $\gamma(t_i)$ by a positive resolution of its double point, put $\epsilon(t_i)=1$.  Otherwise $\epsilon(t_i)=-1$.  

For any generic $\gamma:[0,1]\rightarrow \mathcal{C}$ define $\delta(\gamma)=\sum_{i=1}^n \epsilon(t_i)\sigma(\gamma(t_i))$.  

The set of singular knots forms the {\it discriminant} $D$ in $\mathcal C.$ The codimension two (with respect to $D$) stratum of the discriminant consists of singular knots with two distinct transverse double points. It is easy to see that
$\delta(\alpha')=0$,
for every small generic loop $\alpha'$ going around the
codimension two stratum. This implies (cf. Arnold~\cite{Arnoldcurves}) 
that if $\gamma$ is a generic loop in $\mathcal C$ that starts
at a nonsingular knot $K$,
then $\delta(\alpha)$
depends only on the element of 
$\pi_1(\mathcal C, K)$ realized by a generic loop $\alpha$.

\begin{lem}  The value of the invariant $\delta$ on $\gamma_\text{fib}$ is $2$.  In particular, no element of the homotopy class $[\gamma_{\text{fib}}]\in \pi_1(\mathcal{C},K)$ is represented by a self-isotopy of $K$, that is, an element of $\pi_1(\mathcal{K},K)$.\label{noisotopy.lem}
\end{lem}
\pp Clearly $\delta(\gamma)=0$ for any $\gamma$ which is homotopic to a loop in $\mathcal{K}$.  However $\delta(\gamma_\text{fib})=2$.  The knot $K=\gamma_\text{fib}(0)$ crosses the fiber over the self-intersection point $p$ of $\nu$ twice.  During the homotopy $\gamma_\text{fib}(t)$, these two points move along the fiber at unit speed in opposite directions because the two loops in $F$ one gets by smoothing $\nu$ at $p$ are orientation reversing.  Therefore $\gamma_\text{fib}(t)$ is singular at two times $t_1$ and $t_2$, and $\delta(\gamma_\text{fib})=2$ because $\epsilon(t_1)$ and $\epsilon(t_2)$ are equal and $\sigma(\gamma_\text{fib}(t_1))=\sigma(\gamma_\text{fib}(t_2))=1$. The last identity holds because the two loops adjacent to a double point of singular knots $\gamma_\text{fib}(t_i),$  $i=1,2,$ project to orientation reversing loops on $F$ and hence are not contractible in $M.$  Hence $\gamma_\text{fib}$ is not homotopic to a loop in $\mathcal{K}$.\qed

\begin{lem}  Let $K'$ be any knot homotopic to $K$ through immersed curves, where $K$ is a knot projecting to  $\nu$.  Let $\beta \in \pi_1(\mathcal{C},K')$.  If $\beta$ is
  represented by a self-isotopy of $K'$, that is, an element of
  $\pi_1(\mathcal{K},K')$, then $\beta$ contains a
  representative of the form $\gamma_{rot}^j\gamma_{\text{fr}}^s$ for some
  integers $j$ and $s$. 
\end{lem}\label{standardform.lem}

\begin{remark}
Note that one can show that $s=0$, but we will not need this stronger version of the lemma. 
\end{remark}

\pp Choose a homotopy $\phi$ from $K'$ to $K$ in $\mathcal{C}$, and write $\beta=\phi\alpha\phi^{-1}$ with $\alpha\in \pi_1(\mathcal{C},K)$. Thus $\delta(\beta)=\delta(\phi)+\delta(\alpha)-\delta(\phi)=\delta(\alpha)$.  By Lemma \ref{example2standardform.lem}, we may write $\alpha=[\gamma_\text{fib}^i\gamma_\text{rot}^j\gamma_\text{fr}^k]$.  Now $\delta(\alpha)=i\delta(\gamma_\text{fib})+j\delta(\gamma_\text{rot})+k\delta(\gamma_\text{fr})$.  By Lemma \ref{noisotopy.lem}, $\delta(\gamma_\text{fib})=2$.  Since $\gamma_\text{rot}$ is homotopic to an isotopy, $\delta(\gamma_\text{rot})=0$.  Lastly, $\delta(\gamma_\text{fr})=0$ since the singular knot which appears during $\gamma_\text{fr}$ can be viewed as an ordered pair of two loops, at least one of which is contractible. For $\beta$ to be homotopic to an isotopy we must have $\delta(\beta)=\delta(\alpha)=0$.  Thus $0=2i+0j+0k$, so $i=0$ while thus far $j$ and $k$ can be any integer. Now $\alpha=[\gamma_\text{rot}^j\gamma_\text{fr}^k]$ and $\beta=\phi[\gamma_\text{rot}^j\gamma_\text{fr}^k]\phi^{-1}$.  But $\phi[\gamma_\text{rot}^j\gamma_\text{fr}^k ]\phi^{-1}$ is homotopic to $[\gamma_\text{rot}^j\gamma_\text{fr}^k]\in \pi_1(\mathcal{FC},K')$, where $\gamma_\text{rot}$ is now viewed as a rotation of the parametrizing circle of $K'$ rather than $K$ and $\gamma_\text{fr}$ is the homotopy which passes through a small kink of $K'$ rather than $K$.  Note that the loop $\gamma_\text{rot}$ is a framed isotopy, while the loop $\gamma_{\text{fr}}$ is not even an isotopy.   

\qed

\begin{theorem}\label{theoremexample1} Let $M$ be the oriented total space of an $S^1$-bundle over a non-orientable surface of genus $g\geq 1$ and let $\nu$ be the solid curve pictured in Figure \ref{example2.fig}.  Let $\mathcal{FC}$ be any homotopy class of framed immersions containing a framed knot $K$ that projects to $\nu.$    For each $k\in \mathbb{Z}$ there exists a nowhere-zero vector field $V_k$ on $M$ such that no knot type in $\mathcal{FC}$ is simple.   In particular, for any $V_k$-transverse knot $K'$ in $\mathcal{FC}$, the $V_k$-transverse knots $K'$ and $K'^k$ are homotopic through $V_k$-transverse immersions, isotopic as framed knots, and {\it not} isotopic through $V_k$-transverse knots.\label{nonorientableex.thm}
\end{theorem}

\pp Let $V_k$ be a vector field on $M$ such that the Euler class
$e_{V_k^\perp}\in H^2(M;\Z)$ of the
$2$-plane bundle $V_k^\perp$ on $M$ is Poincar\'e dual to the class
$2k[\tilde{d}]\in H_1(M;\Z)$ where $\tilde{d}$ is some curve in $M$ projecting to
the dotted curve in Figure \ref{example2.fig}.  The fact that such a vector field
exists follows, for example, from the Pontryagin-Thom Construction, see for example Section~\ref{vectorfields.sec}.    Furthermore, choose
$\tilde{d}$ so that it intersects the fiber over $p=\Im \nu \cap \Im
d$ transversely in one point, and is disjoint from $K$. 

To simplify notation in the remainder of the proof, we write $V$ rather than $V_k$. 

 We begin by showing that the knot type of $K$ is simple.  

 We claim that
$k\in \Im \bar{h}_{V}$.  To show this, we must find a framed self-homotopy
$\alpha\in\pi_1(\mathcal{FC}, K)$ such that
$e_{V^\perp}(\alpha_*[S^1\times S^1]])=2k$.  Let
$\alpha=\gamma_{\text{fib}}$.  The lift $\tilde{d}$ of $d$ intersects
$\Im \gamma_{\text{fib}}$ transversely in one point.  Since
$e_{V^\perp}$ is Poincar\'e dual to $2k[\tilde{d}]$, we have
$e_{V^\perp}(\alpha_*[S^1\times S^1]])=2k$.

Now we must show that $k\notin  \Im \bar{i}_{V}$.  It suffices to show $k\notin \Im i_{V}$.  Suppose $\alpha\in \pi_1(\mathcal{K},K)$ and $i_{V}(\alpha)=k$.  We may apply Lemma~\ref{standardform.lem}, and use the lemma in the case where $K'=K$. We conclude that $\alpha$ is represented by the loop $\gamma_{\text{rot}}^j\gamma_{\text{fr}}^s$ for some integers $j$ and $s$.  But $e_{V^\perp}([\gamma_{\text{rot}}^j\gamma_{\text{fr}}^s])=0$ since the corresponding torus is disjoint from $\tilde{d}$, contradicting our assumption.

Now we will show that any knot
type in $\mathcal{FC}$ is simple. Let $K'$ be any knot homotopic to $K$ through framed immersions.  We first check that $k\in \Im \bar{h}_{V}$. (Now the homomorphisms are viewed as maps using basepoint $K'$, and the homomorphism $\bar{i}_V$ is a map from $\pi_1(\mathcal {FK'},K')\rightarrow \Z$, where $\mathcal{ FK'}$ is the connected component of the space of framed knots containing $K'.$)
 Again let $\alpha=\gamma_{\text{fib},K'}$, where now $\gamma_{\text{fib},K'}$ is the framed self-homotopy of $K'$ given by sliding $K'$ around the fiber.  Let $\phi$ be any path from $K'$ to $K$ in $\mathcal{FC}$.  We can write $\alpha=\phi \gamma_{\text{rot},K}\phi^{-1}$, where $\gamma_{\text{rot},K}$ is the self-homotopy of $K$ given by sliding $K$ around the fiber.  Now $e_{V^\perp}(\phi \gamma_{\text{rot},K}\phi^{-1})=2k$ as desired.

Last we show $k\notin \Im \bar{i}_{V}$, and as before it suffices to show that $k\notin \Im i_{V}$.  Suppose $\beta\in \pi_1(\mathcal{K}',K')$, where $\mathcal{K}'$ is the isotopy class of $K'$, and $i_{V}(\beta)=k$.  We may apply Lemma~\ref{standardform.lem}, and conclude that $\beta$ contains a representative of the from $\gamma_{\text{rot},K'}^j\gamma_{\text{fr}}^s$.  Again this is disjoint from $\tilde{d}$ so $e_{V^\perp}(\beta)=0$, contradicting our assumpton.
\qed

The knots in Theorem~\ref{theoremexample1} can be chosen to be Legendrian with an overtwisted complement. This is because the result of Eliashberg~\cite{Eliashbergovertwisted} says that every $2$-plane distribution in $TM$ is homotopic to an overtwisted contact structure.  We make this precise in the following corollary.

\begin{cor} \label{example1cor} Let $(M,\xi_k)$ be an oriented 3-manifold as in Theorem ~\ref{theoremexample1} with a cooriented contact structure $\xi_k$ homotopic to the 2-plane field $V_k^\perp$.   For each framed knot type $\mathcal{FK}\subset\mathcal{FC}$ containing a Legendrian representative $L$, the stabilized Legendrian knots $L_k$ and $L_{-k}$ are isotopic as framed knots, and homotopic as Legendrian immersions, but not isotopic as Legendrian knots; these Legendrian knots can be chosen to be loose.  
	\end{cor}

\begin{proof} Let $L$ be a Legendrian curve which projects to the knot $\nu$ in ~\ref{theoremexample1}.  It is straightforward to check that $L_k$ and $L_{-k}$ are isotopic as framed knots. Recall from the introduction that $L_{k}$ is isotopic to $(L_{-k})^k$ through $V_k$-transverse (in this case, pseudo-Legendrian)
 knots. By Theorem ~\ref{theoremexample1}, $L_{-k}$ and $(L_{-k})^k$ are homotopic through $V_k$-transverse immersions, but not isotopic through $V_k$-transverse knots.   Hence the same is true of $L_{-k}$ and $L_{k}$.  By the $h$-principle for Legendrian immersions~\cite{Gromov} (see Theorem ~\ref{hprlegimm.thm}), $L_{-k}$ and $L_k$ are homotopic as Legendrian immersed curves.  The corollary follows. \end{proof}

\section{Homotopy classes containing Simple and Nonsimple Knot Types}\label{simpleandnonsimple.sec}

In this section we construct framed knot types $\mathcal{FK}_1$ and $\mathcal{FK}_2$ with homotopic representatives ($\mathcal{FK}_i \subset \mathcal{FC}$ for some component $\mathcal{FC}$ of the space of framed curves), such that $\mathcal{FK}_1$ is simple while $\mathcal{FK}_2$ is not.  This is in contrast to the example in the last section, where we described a homotopy class in which all knot types were nonsimple.

\begin{figure}[htbp]
\includegraphics[width=4in]{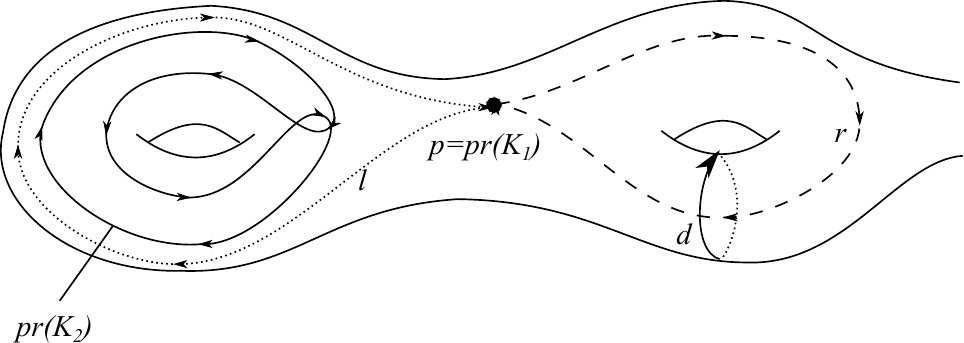}
\caption{$M$ is an $S^1$-bundle over $F$ with vector field $V_k$ such that $e_{V_k^\perp}$ is Poincar\'e dual to $2k[\tilde{d}]\in H_1(M;\Z)$.}\label{simpleandnonsimple.fig}
\end{figure}

Let $M$ be an $S^1$-bundle over an oriented surface of genus $g\geq 2$.  Let $\tilde{d}$ be any curve in $M$ projecting to the solid loop $d$ pictured in Figure \ref{simpleandnonsimple.fig} to $M$ and let $V_k$ be a vector field on $M$ such that $e_{V_k^\perp}$ is Poincar\'e dual to $2k[\tilde{d}]\in H_1(M;\Z)$. Let $K_1$ be a vertical fiber with some framing.

\begin{figure}[htbp]
\includegraphics[width=2in]{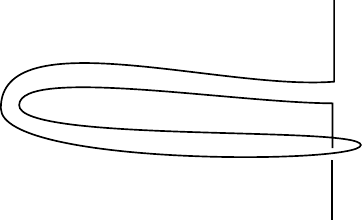}
\caption{$K_2$ is obtained from $K_1$ by pulling a small arc of $K_1$ around a loop in $M$ which projects to a loop freely homotopic to $l^{-1}$ on $F$.}\label{K2.fig}
\end{figure}

Let $K_2$ be a framed knot obtained from $K_1$ by pulling a small arc of $K_1$ around a loop in $M$ which projects to a loop freely homotopic to $l^{-1}$ on $F$ and then passing through a double point,  see Figure \ref{K2.fig}.  The loop $l$ is pictured in Figures \ref{simpleandnonsimple.fig}.
Let $\mathcal{C}$ denote the component of the space of immersions containing $K_1$ and $K_2$.  Now we characterize all possible self-homotopies of $K_1$.  Recall that the definition of $\gamma_\rho$ can be found in Subsection \ref{gammarho}.
\begin{lem}[\cite{ChernovFramed} proof of Lemma 6.11]
Let $\alpha\in \pi_1(\mathcal{C},K_1)$.  Then $\alpha$ contains a representative of the form $\gamma_{\rho}\gamma_\text{rot}^j\gamma_\text{fr}^k$ for some choice of $\gamma_\rho$, where $\rho\in \pi_1(F,p)$.\label{standardformexample3.lem}
\end{lem}
\pp The proof is similar to that of Lemma \ref{example2standardform.lem}.  Again the trace of the basepoint under a self homotopy commutes with the curve $K_1$.  In this case, the trace of the basepoint of $\alpha$ must commute with the class of the fiber $[f]\in \pi_1(STM, \vec{K}(1))=\pi_1(M,K(1))$.  But $[f]$ is in the center of $\pi_1(M, K(1))$, so $\tau(\alpha)$ can be any element of $\pi_1(M, K(1))$.  In particular we may write $\tau(\alpha)=\tilde{\rho}\cdot[f]^j$ for some loop $\tilde{\rho}$ which projects to $\rho\in \pi_1(F,p)$.  But for some choice of loop $\gamma_\rho$, the loop $[\gamma_\rho\gamma_\text{rot}^j]$ has trace $\tilde{\rho}\cdot [f]^j$ as well.  Now the lemma follows from Proposition \ref{sametrace.prop}.\qed

\begin{prop}  Let $\phi:[0,1]\rightarrow \mathcal{C}$ be a path from $K_2$ to $K_1$ which unclasps $K_2$ and is an isotopy at all other times. Every loop in $\pi_1(\mathcal{C},K_2)$ has a representative of the form $\phi \gamma_{\rho}\gamma_\text{rot}^j\gamma_\text{fr}^k \phi^{-1}$ for some $\rho \in \pi_1(F,p)$.\label{K2stdform}
\end{prop}

\pp This follows directly from Lemma \ref{standardformexample3.lem}.\qed

\begin{defin}[Schneiderman's invariant]
Now we recall an invariant due to Schneiderman~\cite{Schneiderman} which we use in the proof of the next lemma.  Let $X$ be a $4$-manifold and let $A: S^1\times [0,1], S^1\times\{0,1\}\rightarrow(X,\partial X)$ be a properly immersed annulus. Let $x$ be a basepoint of $X$ and let $a$ be a basepoint of $\Im A$. A {\it whisker} for $A$ is a choice of path from $x$ to $a$; fix some whisker $\omega$.  For each self-intersection point $p$ of $A$, the {\it sheets} at $p$ are the two transversely intersecting  immersed $2$-disks in a small neighborhood of $p$ in $A$.  For each self-intersection $p$ define a loop $g_p\in \pi_1(X,x)$ as follows: go along $\omega$ to $a$, go along a path in $A$ to $p$, switch sheets, return to $a$ without passing through any other double points of $A$, and then return to $x$ along $\omega^{-1}$. The loop $g_p$ is well defined up to powers the loop $\kappa=\omega A_*(S^1\times t_a)\omega^{-1}$ where $a\in A(S^1\times t_a)$.  Define a sign $\epsilon(p)$ by comparing the orientation of $X$ at $p$ with the orientation given by the two sheets of $A$ at $p$.  Now 
$$\mu(A)=\sum_{p\in A\cap A}\epsilon(p) [g_p].$$
Let $\Lambda_\kappa=\Z[\pi_1(X,x)]/\{g-\kappa^n g^{\pm 1} \kappa^m\}$ where $\Z[\pi_1(X,x)]$ denotes the free abelian group generated by the elements of $\pi_1(X,x)$.  Note that if one wants an invariant of homotopy rather than just regular homotopy one should add $\Z[1]$ to the denominator of the quotient; for our purposes a regular homotopy invariant is enough.

Following Wall ~\cite{Wall}, Schneiderman \cite[Proposition 4.1.2]{Schneiderman} proves that $\mu(A)$, when viewed as an element of the quotient $\Lambda_\kappa$, is an invariant of regular homotopy, and whenever $\mu$ vanishes on $A$, the double points of $A$ can be paired off with Whitney disks.  In higher dimensions, because of the Whitney trick, $\mu$ vanishes if and only if  $A$ is regularly homotopic to an embedding; in dimension $4$, $\mu$ vanishing is just a necessary condition for $A$ to be regularly homotopic to an embedding.
\begin{lem} The class $[\phi \gamma_{\rho}\gamma_\text{rot}^j\gamma_\text{fr}^k \phi^{-1}]\in \pi_1(\mathcal{C},K_2)$ is not represented by an element of $\pi_1(\mathcal{K}_2,K_2)$ if $\rho\neq[l]^s$ for some $s\in \Z$.\label{immersedannulus.lem}
\end{lem}
\pp  Let $x_1$ and $x_2$ be the preimages of the double point of the singular knot which appears during the homotopy $\phi$ at time $s_0$; that is, $\phi_{s_0}(x_1)=\phi_{s_0}(x_2)$ for some $s_0\in[0,1]$.  We assume that this unclasping homotopy has the properties that $\text{pr}(\phi_s(x_1)\equiv p)$, and $\text{pr}(\phi_s(x_2))=l(s)$ for $s\in[s_0,1]$. The loop $l$ is pictured in Figure \ref{simpleandnonsimple.fig}. 
\end{defin}

By modifying our choice of loop $\gamma_\rho$, we may assume $j=0$.  Also, $\gamma_\text{fr}$ commutes with $\gamma^\rho$ and $\phi$, so we begin by considering $\phi\gamma_\rho\phi^{-1}$. 

We consider the track of the homotopy $\phi\gamma_\rho \phi^{-1}$ in $[0,1]\times M$.  By this we mean, consider the following map  $A:[0,1]\times S^1\rightarrow [0,1]\times M$:
$$A(s,t)=
\left\{\begin{array}{ll}s\times \phi_{3s}(t)&\text{for }s\in [0,1/3]\\ s\times (\gamma_\rho)_{3s-1}(t)&\text{for } s\in[1/3,2/3]\\
s\times \phi^{-1}_{3s-2}(t)&\text{for } s\in [2/3,1] \end{array}\right.$$

\begin{figure}[htbp]\includegraphics[width=13cm]{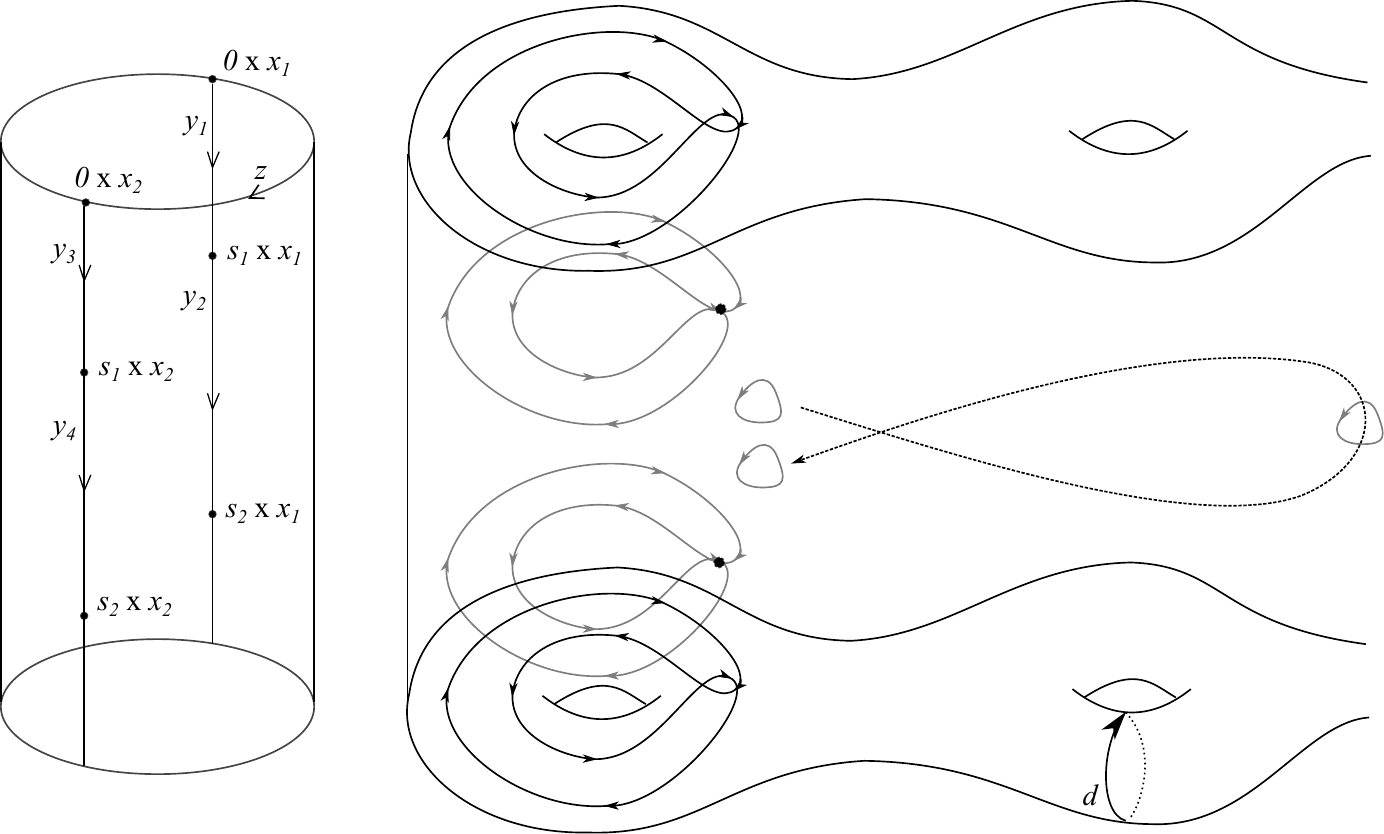}
\caption{The track of the homotopy $\phi\gamma_\rho\phi^{-1}$ in $[0,1]\times M$. }\label{annulus.fig}
\end{figure}

  If the loop $\phi\gamma_\rho \phi^{-1}$ is homotopic to a loop in $\pi_1(\mathcal{K}_2,K_2)$, then the immersed annulus $A:[0,1]\times S^1\rightarrow [0,1]\times M$ is homotopic to an embedded one, through maps fixing the boundary.  We compute $\mu(A)$ and conclude that if $\rho$ does not commute with $[l]$ in $\pi_1(F,p)$, then $A$ is not homotopic to an embedded annulus.

Now we compute the two terms of $\mu(A)$ corresponding to the times $s_1$ and $s_2$ at which the self-homotopy $\phi \gamma_\rho \phi^{-1}$ of $K_2$ is singular.  For our purposes it will actually be enough to compute the projections of the terms of $\mu(A)$ to $F$.

The projection to $F$ of the term of $\mu (A)$ corresponding to the first singular time $s_1$ is
$$\text{pr}_*(A(y_1)A(y_3^{-1}z^{-1}))=[l].$$

\begin{figure}\includegraphics[width=2in]{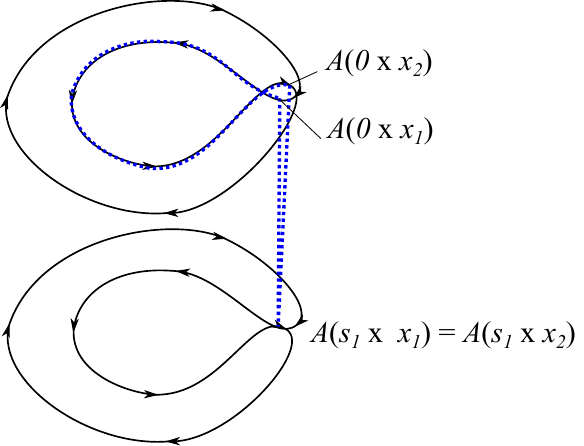}
\caption{The term $g_{p\times t_1}$ projects to $\beta^{-1}$ on $F$.}
\label{Firstterm.fig}
\end{figure}

A picture of this loop in $[0,1]\times M$ can be found in Figure \ref{Firstterm.fig}.
Now we consider the projection to $F$ of the term of $\mu(A)$ corresponding to the second singular time $s_2$.

The projection to $F$ of the term of $\mu(A)$ corresponding to the second singular time $s_2$ is
$$\text{pr}_*(A(y_1y_2)A(y_4^{-1}y_3^{-1}z^{-1})).$$

Since $\text{pr}(\phi_s(x_1)\equiv p)$, and $\text{pr}(\phi_s(x_2))=l(s)$ for $s\in[s_0,1]$, we have
$$\text{pr}_*(A(y_2))=\rho$$
and
$$\text{pr}_*(A(y_4))=[l]\rho [l^{-1}].$$

Thus 
$$\text{pr}_*(A(y_1y_2)A(y_4^{-1}y_3^{-1}z^{-1}))=\rho[l]\rho^{-1}[l]^{-1}[l].$$

Now we check whether the two terms $[l]$ and $\rho[l]\rho^{-1}[l]^{-1}[l]$ are equal in the quotient $\text{pr}_*\Lambda_\kappa$, where $\kappa=A_*(0\times S^1)=K_2$.  Since $K_2$ is homotopic to the fiber, $\text{pr}_*\Lambda_\kappa=\Z[\pi_1(F,p)]/\{g-g^{\pm 1}\}.$  It follows that these two terms of $\text{pr}_*(\mu(A))$ do not cancel unless $[l]$ and $\rho$ commute in $\pi_1(F,p)$.  Since centralizers in the fundamental group are infinite cyclic, and since $[l]$ is not a nontrival power of another loop, $[l]$ and $\rho$ commute if and only if $\rho=[l]^i$ for some $i\in \Z$.\qed
\begin{theorem}\label{Theoremexample2}
 Let $M$ be an $S^1$-bundle over an orientable surface of genus $g\geq 2$.  Let $\mathcal{FC}$ be the framed homotopy class of the $S^1$ fiber with any framing.  Then for any $k\in \Z$ there exists a nowhere-zero vector field $V_k$ such that the knot type of $K_1$ is simple in $(M, V_k)$ while the knot type of $K_2$ is not.  In particular, the knots $K_2$ and $K_2^k$  are homotopic through $V_k$-transverse immersions, isotopic as framed knots, and {\it not} isotopic through $V_k$-transverse knots.\label{orientableex.thm}
\end{theorem}

\pp  In this proof, we write $V$ rather than $V_k$ to simplify notation.

First we consider the knot type of $K_1$.  The loop $\gamma_{[r]}$ is a framed self-homotopy and a framed self-isotopy of $K_1$.  The corresponding torus intersects $[\tilde{d}]$ once transversely.  Therefore $e_{V^\perp}(\gamma_{[r]*}[S^1\times S^1])=2k$, so $k\in \Im \bar{h}_V$ and $k\in \Im \bar{i}_V$.  Furthermore $k$ is the smallest positive integer in $\Im \bar{h}_V$ or $\Im \bar{i}_V$, since any such torus intersects $[\tilde{d}]$ at least once or not at all.  This $\Im \bar{h}_V=\Im \bar{i}_V=k\Z$ and the knot type $K_1$ is simple by Corollary \ref{equalimages.cor}.

Now we consider the knot type of $K_2$.  In this case the loop $\phi\gamma_{[r]}\phi^{-1}$ is a framed self-homotopy of $K_2$ so as before $\Im \bar{h}_V=k\Z$.  Let $\alpha$ be an arbitrary framed self-isotopy of $K_2$.  By Proposition \ref{K2stdform} we know $\alpha$ is represented by a loop of the form $\phi \gamma_{\rho}\gamma_\text{rot}^j\gamma_\text{fr}^k \phi^{-1}$ for some $\rho \in \pi_1(F,p)$.  By Lemma \ref{immersedannulus.lem}, since $\alpha$ is a framed self-isotopy, $\rho=[l]^i$ for some $i\in \Z$.  But $\phi \gamma_{[l]^i}\gamma_\text{rot}^j\gamma_\text{fr}^k \phi^{-1}$ is disjoint from $[\tilde{d}]$, so $e_{V^\perp}(\alpha_*[S^1\times S^1])= 0$.  Hence $\Im i_V=0$, and the knot type of $K_2$ is not simple.  In particular Lemmas \ref{Vtranshomo.lem} and \ref{Vtransiso.lem} imply $K_2$ and $K_2^k$ are homotopic through $V$-transverse immersions and not isotopic as $V$-transverse knots.
\qed

As in Theorem~\ref{theoremexample1}, the knots in Theorem~\ref{Theoremexample2} can be chosen to be Legendrian with an overtwisted complement. The proof of Corollary \ref{example2cor} is similar to the proof of Corollary \ref{example1cor}.
\begin{cor} \label{example2cor} Let $(M,\xi_k)$ be an oriented 3-manifold as in Theorem ~\ref{Theoremexample2}, with a cooriented contact structure $\xi_k$ homotopic to the 2-plane field $V_k^\perp$. Let $L$ be a Legendrian knot which is smoothly isotopic to the knot $K_2$ in Theorem ~\ref{Theoremexample2}. Then the stabilized Legendrian knots $L_k$ and $L_{-k}$ are isotopic as framed knots, and homotopic as Legendrian immersions, but not isotopic as Legendrian knots; these Legendrian knots can be chosen to be loose.  
\end{cor}

\section{Applications to Legendrian Knot Theory}\label{legendrian.sec}

A Legendrian knot $L$ in an overtwisted contact $3$-manifold $(M,\xi)$ is called {\it loose} if $L\subset M\setminus D$ where $D$ is an overtwisted disk.  In many situations, loose Legendrian knots are completely determined by classical invariants.  Dymara \cite[Theorem 4.1]{Dymara} showed this is true, for example, when $\xi$ is trivializable and $|L|=\infty$.  Ding and Geiges \cite{DingGeiges} generalized Dymara's result to the case where the connected component of the space of immersions containing $L$ contains infinitely many components of the space of Legendrian knots, and again $|L|=\infty$.  

Our theory leads to a classification of loose Legendrian knots in terms of generalized classical invariants which generalizes the results of Dymara and Ding-Geiges.  We use an $h$-principle stated by Cieliebak and Eliashberg, and attributed to Dymara \cite{Dymara} and Eliashberg-Fraser \cite{EliashbergFraser}. 

\subsection*{$h$-principles for Legendrian knots}  We recall some terminology from Cieliebak and Eliashberg \cite{CieliebakEliashberg}, using their notation.  If $M$ and $N$ are manifolds, a {\it monomorphism} is a bundle map $F:TM\rightarrow TN$ which is injective on each fiber.  For example, given an immersion $f:S^1\rightarrow N$, the differential $df$ is a monomorphism $TS^1\rightarrow TN$.

Now let $(M,\xi)$ be a contact $2n+1$-manifold and $\Lambda$ a manifold of dimension $\leq n$.  In our case $n=1$ and $\Lambda=S^1$.    A {\it formal Legendrian embedding} is a pair $(f,F^s)$ consisting of a smooth embedding $f:\Lambda\rightarrow M$ and a homotopy of monomorphisms $F^s$ over $f$ starting at $F^0=df$ and ending at an {\it isotropic} or {\it Legendrian} monomorphism, meaning $F^1$ lies in $\xi$.

Two formal Legendrian embeddings are called {\it formally isotopic} if they are isotopic as (connected by a path of) formal Legendrian embeddings.

\begin{thm}[\cite{CieliebakEliashberg}Theorem 7.19 b]  Let $(M,\xi)$ be a closed connected overtwisted contact $3$-manifold and $D\subset M$ an overtwisted disk.

Let $(L_t,F^s_t)$, $s, t\in[0,1]$ be a formal Legendrian isotopy in $M$ connecting two genuine Legendrian embeddings $L_0,L_1:S^1 \rightarrow M\setminus D$.  Then there exists a Legendrian isotopy $\tilde{L}_t:S^1\rightarrow M\setminus D$ connecting $\tilde{L}_0=L_0$ and $\tilde{L}_1=L_1$, which is homotopic to $(L_t,F^s_t)$ through formal Legendrian isotopies with fixed endpoints.\label{cethm}
\end{thm}

We use Theorem \ref{cethm} to obtain the following result.

\begin{thm} Let $(M,\xi)$ be a closed  overtwisted contact manifold with overtwisted disk $D$, and let $V$ be the coorienting vector field of $\xi$.  Let $L_1$ and $L_2$ be Legendrian knots in $M\setminus D$ that are $V$-transverse isotopic.  Then they are Legendrian isotopic in $M\setminus D$.\label{Vtransimpliesformal}
\end{thm}

\begin{proof}  Fix an auxiliary Riemannian metric with the property that $V$ is everywhere orthogonal to $\xi$.  Let $L^V_t$ denote the $V$-transverse isotopy from $L_1$ to $L_2$.  That is, for all $t\in[0,1]$, $L^V_t$ is a $V$-transverse knot, and $L^V_0=L_1$ and $L^V_1=L_2$ are Legendrian. We use the $V$-transverse isotopy $L^V_t$ from $L_1$ to $L_2$ to construct a formal Legendrian isotopy from $L_1$ to $L_2$.  Let $P_u=\text{proj}_{\xi}(L^V_t)'(u)$ denote the normalized orthogonal projection of the tangent vector $(L^V_t)'(u)$ to the contact plane $\xi_{L^V_t(u)}$.  Let $V_u$ denote the vector in $V$ at the point $L^V_t(u)$.  Let $\alpha_u$ denote the angle between $(L^V_t)'(u)$ and its projection $P_u$.  Let
$$V^{s,t}(u)=\cos((1-s)\alpha-u)P_u+\sin((1-s)\alpha_u)V_u,$$
which is a unit vector in $T_{L^V_t(u)}M$.  Then $V^{0,t}(u)$ is equal to $(L^V_t)'(u)$ and $V^{1,t}$ is tangent to $\xi_{L^V_t(u)}$.  Hence for each $t$, the pair $(L^V_t(u),V^{s,t}(u))$ is a formal Legendrian embedding.  Now the theorem follows from Theorem \ref{cethm}.
\end{proof}

As an application of the above theorem we get the following:
\begin{thm} Let $(M,\xi)$ be a closed connected overtwisted contact manifold with overtwisted disk $D$.  Let $L_1, L_2\subset M\setminus D$ be two loose Legendrian knots in the smooth isotopy class $\mathcal{K}$.  Assume the following conditions hold:

\begin{enumerate}
\item $L_1$ and $L_2$ are isotopic as framed knots (and hence lie in some framed isotopy class $\mathcal{FK}$)
\item $L_1$ and $L_2$ are homotopic as Legendrian immersions
\item $\Im \bar{h}_V=\Im \bar{i}_V$, where $V$ is a coorienting vector field for $\xi$, and $\bar{h}_V$ and $\bar{i}_V$ are the kink-cancelling homomorphisms associated to $\mathcal{FK}$.
\end{enumerate}
Then $L_1$ and $L_2$ are isotopic as Legendrian knots.\label{loosethm}
\end{thm}
\begin{proof} Because $L_1$ and $L_2$ are Legendrian homotopic they must be $V$-transverse homotopy.  Since $L_1$ and $L_2$ are also in the same framed isotopy class, Corollary \ref{equalimages.cor} implies $L_1$ and $L_2$ are $V$-transverse isotopic.  Finally, Theorem \ref{Vtransimpliesformal} implies $L_1$ and $L_2$ are Legendrian isotopic.
\end{proof}
Below we state theorems of Dymara \cite{Dymara} and Ding-Geiges \cite{DingGeiges} which we will show follow from Theorem \ref{loosethm}.  

\begin{thm}[\cite{Dymara}, Theorem 4.1] Let $(M,\xi)$ be a contact manifold with an overtwisted disk $D$ and a trivializable contact bundle.  Let $L_1, L_2\subset M\setminus D$ be two loose Legendrian knots in the smooth isotopy class $\mathcal{K}$.  Assume the following conditions hold:
\begin{enumerate} 
\item $|K|=\infty$ 
\item $L_1$ and $L_2$ are isotopic as framed knots
\item the rotation numbers of $L_1$ and $L_2$ with respect to some trivialization of $\xi$ are equal.
\end{enumerate}
Then $L_1$ and $L_2$ are isotopic as Legendrian knots.\label{dymaraloosethm} 
\end{thm}

This theorem of Dymara in the case of closed $M$ can be viewed as a special case of Theorem \ref{loosethm}.  This is because when $\xi$ is trivializable, $\bar{h}_V=\bar{i}_V=0$ for every framed isotopy class.  Note that for trivializable $\xi$ the $V$-transverse homotopy classes of curves in a fixed component of the space of framed curves $\mathcal{FC}$ are enumerated by a $\mathbb{Z}$-valued rotation number, which is obtained by projecting the velocity vector of a $V$-transverse curve to the planes of $\xi$.

Ding and Geiges \cite[Theorem 6]{DingGeiges} generalized Dymara's theorem to the case where $\xi$ is not necessarily trivializable:

\begin{thm}[\cite{DingGeiges}, Theorem 6]  Let $(M,\xi)$ be a closed connected contact manifold with an overtwisted disk $D$.  Let $L_1, L_2\subset M\setminus D$ be two loose Legendrian knots in the smooth isotopy class $\mathcal{K}$.  Assume the following conditions hold:
\begin{enumerate}
\item $|K|=\infty$
\item $L_1$ and $L_2$ are isotopic as framed knots
\item $L_1$ and $L_2$ are homotopic as Legendrian immersions
\item the connected component of the space of framed immersions containing $L_1$ and $L_2$ contains infinitely many components of the space of Legendrian immersions.
\end{enumerate}
Then $L_1$ and $L_2$ are isotopic as Legendrian knots.
\end{thm}

As shown in \cite[Proposition 3.1.4]{TchernovVassiliev}, the assumption that the connected component of the space of framed immersions $\mathcal{CF}$ contains infinitely many components of the space of Legendrian immersions is equivalent to the following: For any $\alpha\in H_2(M;\mathbb{Z})$ which is realizable by a mapping $\mu:S^1\times S^1\rightarrow M$ with meridian freely homotopic to a loop in $\mathcal{FC}$, we have $e_{\xi}(\alpha)=0$.  Therefore assumption $(4)$ of Ding and Geiges implies that $\bar{h}_V=\bar{i}_V=0$.  Thus their theorem is also a special case of Theorem \ref{loosethm}.

Ding and Geiges \cite[Remark, p. 121]{DingGeiges} remark that conditions $(1)$ and $(4)$ of their theorem are necessary unless one makes certain ad-hoc assumptions about $L_1$ and $L_2$.

It is possible to coarsely classify zero-homologous loose knots even when two such knots do not have the same overtwisted disk in their complement. The following classification of loose knots up to contactomorphism was given in the work of Etnyre~\cite[Theorem 1.4]{Etnyreovertwisted}. (According to~\cite{Etnyreovertwisted} different proofs of this result were independently obtained by Geiges and Klukas.)
Zero-homologous framed knots corresponding to a given unframed knot $K$ are enumerated by the self-linking number, which is the Thurston-Bennequin invariant $\tb$ of a Legendrian knot with the natural framing; and the sum of the 
$\tb$ and $\rot$ of a Legendrian knot is always odd. 

\begin{thm}[\cite{Etnyreovertwisted}, Theorem 1.4]
Let $(M,\xi)$ be an overtwisted contact manifold. For each zero-homologous knot type $\mathcal K$ and a pair of integers $(t,s)$ such that $t+s$ is odd, there is a unique, up to contactomorphism, loose Legendrian knot in $\mathcal K$ satsifying $\tb (K)=t$ and $\rot(K)=r.$
\end{thm}

\section{Appendix: $h$-principles}\label{hprinciples.sec}

\subsection*{$h$-principle for immersed curves}
We will use an $h$-principle in order to understand the topology of a component $\mathcal{C}$ of the space of framed immersions of the circle into $M$.  Let $p:STM\rightarrow M$ be the unit two-sphere tangent bundle over $M$.  

\begin{thm}[Hirsch-Smale $h$-principle, \cite{Gromov}; see also Theorem 7.1 of \cite{CieliebakEliashberg}]  Let $M$ be a $3$-manifold. The space of immersed curves in $M$ is weak homotopy equivalent to the space $\Omega STM$ of continuous free loops in $STM$.  The weak homotopy equivalence is given by mapping the immersed curve $C:S^1\rightarrow M$ to the loop $\vec{C}\in \Omega STM$, where $\vec{C}: t\mapsto C'(t)$, $t\in S^1$.
\end{thm}

Let $a:S^1\rightarrow \Omega STM$, with $a(1)(t)=\omega(t):S^1\rightarrow STM$. Then $a$ can also be viewed as a map $a:S^1\times S^1\rightarrow STM$ defined by $a(s\times t)=a(s)(t)$.  Identify $S^1$ with the unit circle in $\mathbb{C}$.  Let $\tau(a)=a(1\times t)$ be the loop in $STM$ traced by the basepoint $\omega(1)$ during the homotopy $a$.

The map $\tau$ is in fact a homomorphism $\tau:\pi_1(\Omega STM, \omega)\rightarrow \pi_1(STM,\omega(1))$. The image of $\tau$ is precisely the centralizer $Z(\omega)<\pi_1(STM,\omega(1))$.  Now we consider its kernel.  Suppose that $\tau([a_1])=\tau([a_2])$.  Then the $1$-skeleta $a_{1*}(1\times S^1)\cup a_{1*}(S^1\times 1)$ and $a_{2*}(1\times S^1)\cup a_{2*}(S^1\times 1)$ are homotopic in $STM$; in fact, they agree on $\omega=a_{i*}(1\times S^1)$.  The obstruction to $a_1$ and $a_2$ being homotopic is therefore the element of $\pi_2(STM)$ formed by gluing the two $2$-cells in the images of $a_1$ and $a_2$ together along their common $1$-skeleta.  Note that $TM$ is trivial, so we choose an identification $STM\simeq S^2\times M$, which gives an isomorphism $\pi_2(STM)\simeq \mathbb{Z}\times \pi_2(M)$. In addition, we can identify $\pi_1(STM,\omega(1))$ with $\pi_1(M,p(\omega(1)))$ and view $\tau$ as a map into $\pi_1(M,p(\omega(1)))$ when convenient.

In our case, $M$ will be an $S^1$-bundle over a surface of genus at least $2$, so $\pi_2(M)=0$.  

\begin{prop}[\cite{ChernovFramed}] Let $M$ be an oriented $3$-manifold with $\pi_2(M)=0$, let $\mathcal{C}$ be a connected component of the space of immersed curves in $M$, and let $K\in \mathcal{C}$ be a knot.  Let $\alpha_1,\alpha_2\in \pi_1(\mathcal{C}, K)$ such that $\tau(\alpha_1)=\tau(\alpha_2)\in \pi_1(M, K(1))$.  Then $\alpha_1\gamma_\text{fr}^m=\alpha_2$ where $m\in \mathbb{Z}$ is the first coordinate in $\pi_2(STM)\simeq \mathbb{Z}\times \pi_2(M)$. 
\label{sametrace.prop}
\end{prop}

\subsection*{h-principle for Legendrian immersions}

\begin{thm}[Gromov \cite{Gromov}]\label{hprlegimm.thm} Let $(M,\xi)$ be an oriented 3-manifold with cooriented contact structure $\xi$, and let $S_\xi M$ denote the corresponding sphere bundle. The space of Legendrian immersions of $S^1$ into $(M,\xi)$ is weak homotopy equivalent to the space $\Omega S_\xi M$ of continuous free loops in $S_\xi M$.  The weak homotopy equivalence is given by mapping the Legendrian immersion $L:S^1\rightarrow M$ to the loop $\vec{L}\in \Omega S_\xi M$, where $\vec{L}:t\mapsto L'(t)$, $t\in S^1$.
\end{thm}

{\bf Acknowledgments.}
This paper was written when the authors were at Max Planck Institute for Mathematics, Bonn and the authors thank the institute for its hospitality. This work was partially supported by a grant from the Simons Foundation
$\#235674$ to Vladimir Chernov and $\#523862$ to Patricia Cahn. 

The authors are thankful to Stefan Nemirovski and Spiridon Adams-Florou for enlightening discussions.
The authors are grateful to the anonymous referee for all the valuable comments and suggestions.


\begin{thebibliography}{99999}
\bibitem{Arnoldcurves}
{\bf V.~I.~Arnold} {\it Plane curves, their invariants, perestroikas and
classifications}, Singularities and Bifurcations (V.I.~Arnold, ed.) Adv.~Sov.~Math.~vol. {\bf 21,} 
(1994), pp.~39--91

\bibitem{BenedettiPetronio1} 
{\bf R.~Benedetti, C.~Petronio} {\it Combed 3-manifolds with concave boundary, framed links, and pseudo-Legendrian links,\/} J. Knot Theory Ramifications {\bf 10} (2001), no. 1, 1--35. 

\bibitem{BenedettiPetronio2} 
{\bf R.~Benedetti, C.~Petronio}
{\it Reidemeister-Turaev torsion of 3-dimensional Euler structures with simple
boundary tangency and pseudo-Legendrian knots,\/}
Manuscripta Math.~{\bf 106} (2001), no. 1, 13--61

\bibitem{CahnChernovSadykov}
{\bf P.~Cahn, V.~Chernov, R.~Sadykov} {\it The number of framings of a knot in a $3$-manifold\/}
preprint arXiv:1404.5851 (2014) 8 pages

\bibitem{ChernovFramed}
{\bf V.~Chernov} {\it Framed knots in 3-manifolds and affine self-linking numbers,} J.~Knot Theory Ramifications {\bf 14} (2005), no. 6, 791--818

\bibitem{ChernovTransverse}
{\bf V.~Chernov} {\it Relative framing of transverse knots,} Int.~Math.~Res.~Not. (2004) no.~52, 2773--2795


\bibitem{CieliebakEliashberg}
{\bf K.~Cieliebak Ya.~Eliashberg} {\it From Stein to Weinstein and Back Symplectic Geometry of Affine Complex Manifolds\/} a book available at
http://www.mathematik.uni-muenchen.de/~kai/research/stein.pdf 372 pages 

\bibitem{DingGeiges}
{\bf F. Ding and H. Geiges} {\it Handle moves in contact surgery diagrams,} J. Topology   {\bf 2},  no.~1,  pp. 105--122 (2009). 

\bibitem{DoCarmo} 
{\bf M.~do Carmo} {\it  Riemannian geometry.\/} Mathematics: Theory \& Applications. Birk\"auser Boston, Inc., Boston, MA (1992)

\bibitem{Dymara}
{\bf K.~Dymara} {\it Legendrian knots in overtwisted contact structures on $S^3$} Ann.~Global Anal.~Geom.~{\bf 19}, no.~3, 293--305 (2001)

\bibitem{DymaraOvertwisted}
{\bf K.~Dymara} {\it Legendrian knots in overtwisted contact structures} http://arxiv.org/abs/math/0410122 (2004) 34 pages

\bibitem{Eliashbergovertwisted}
{\bf Ya.~Eliashberg} {\it Classification of overtwisted contact structures on 3-manifolds} Invent.~Math.~{\bf 98} (1989), no. 3, 623--637


\bibitem{Eliashberg}
{\bf Ya.~Eliashberg} {\it Contact $3$-manifolds twenty years since J.~Martinet's
work,\/} Ann. Inst. Fourier (Grenoble) {\bf 42} no. 1-2 (1992), pp. 165--192

\bibitem{EliashbergFraser}
{\bf Ya.~Eliashberg and M.~Fraser} {\it Topologically trivial Legendrian knots} J.~Symp.~Geom. {\bf 7}
no. 2, 77--127 (2009).


\bibitem{Etnyreovertwisted}
{\bf J.~Etnyre}
{\it On knots in overtwisted contact structures} 
Quantum Topol.~{\bf 4} (2013), no.~3, 229--264

\bibitem{Geiges}
{\bf H.~Geiges} {\it An introduction to contact topology.\/} Cambridge Studies in Advanced Mathematics, {\bf 109} Cambridge University Press, Cambridge (2008) 

\bibitem{Gromov}
{\bf M. Gromov} {\it Partial Differential Relations,} Springer, Berlin, Heidelberg, 1986.


\bibitem{Lutz}
{\bf R.~Lutz} {\it Structures de contact sur les fibr\'es principaux en cercles
de dimension $3$,\/} Ann.~Inst.~Fourier, {\bf 3} (1977), pp. 1--15

\bibitem{McCullough}
{\bf D. McCullough} {\it Homeomorphisms which are Dehn twists on the boundary,} Algebr. Geom. Topol., Vol {\bf 6} (2006), pp. 1331--1340.


\bibitem{Schneiderman}
{\bf R. Schneiderman.} {\it Algebraic linking numbers of knots in 3-manifolds,} Algebr. Geom. Topol. {\bf 3} (2003), pp. 921--968.

\bibitem{TchernovpseudoLegendrian}
{\bf V.~Tchernov (Chernov).} {\it Isomorphism of the groups of Vassiliev invariants of Legendrian and pseudo-Legendrian knots,} Compositio Math.~{\bf 135} (2003), no. 1, 103--122

\bibitem{TchernovVassiliev}
{\bf V.~Tchernov (Chernov).}  {\it Vassiliev invariants of Legendrian, transverse, and framed knots in contact three-manifolds,} Topology {\bf 42} 
(2003), no. 1, 1--33

\bibitem{Trace}
{\bf B. Trace.} {\it On the Reidemeister Moves of a Classical Knot,}  Proc. Amer. Math. Soc. {\bf 89} (1983), no. 4, 722--724.   

\bibitem{Wall}
{\bf C. T. C. Wall.} {\it Surgery on Compact Manifolds.} London Math.~Soc.~Mono-
graphs 1, Academic Press, (1970) or Second Edition, edited by A. Ranicki,
Math. Surveys and Monographs {\bf 69}, A.M.S.

\end{thebibliography}
\end{document}